\documentclass[twocolumn]{autart}    

\usepackage{graphicx}          
\usepackage{cite}
\usepackage{amsmath,amssymb,amsfonts}

\usepackage{algorithm,algorithmic}
\usepackage{hyperref}
\usepackage{textcomp}
\usepackage{epstopdf}
\usepackage{stfloats}
\usepackage{subfigure}       
\usepackage{xcolor}    
\usepackage{chngcntr}   
\usepackage{makecell}
\usepackage{tabularx} 
\usepackage{multirow}

\newtheorem{assumption}{Assumption}

\begin{document}
	
	\begin{frontmatter}
		
		\title{Distributed Constraint-coupled Resource Allocation: Anytime Feasibility and Violation Robustness} 
		
		\thanks[footnoteinfo]{This work was supported in part by National Key R\&D Program of China under the grant 2022YFB3303900, and the NSF of China under the Grants 62173225. Corresponding author Shanying Zhu.}
		
		\author[sjtu1,sjtu2,sjtu3]{Wenwen Wu}\ead{wuwenwen@sjtu.edu.cn},   
		\author[sjtu1,sjtu2,sjtu3]{Shanying Zhu}\ead{shyzhu@sjtu.edu.cn},
		\author[sjtu1,sjtu2,sjtu3]{Cailian Chen}\ead{cailianchen@sjtu.edu.cn},       
		\author[sjtu1,sjtu2,sjtu3]{Xinping Guan}\ead{xpguan@sjtu.edu.cn}
		
		\address[sjtu1]{School of Automation and Intelligent Sensing, Shanghai Jiao Tong University, Shanghai 200240, China;}
		\address[sjtu2]{Key Laboratory of System Control and Information Processing, Ministry of Education of China, Shanghai 200240, China;}
		\address[sjtu3]{Shanghai Key Laboratory of Perception and Control in Industrial Network Systems, Shanghai, 200240, China;}

		\begin{keyword}                           
			Distributed optimization; real-time system; control barrier function; primal-dual method.               				  
		\end{keyword}                             

		\begin{abstract}                          
			This paper considers distributed resource allocation problems (DRAPs) with a coupled constraint for real-time systems. Based on primal-dual methods, we adopt a control perspective for optimization algorithm design by synthesizing a safe feedback controller using control barrier functions to enforce constraint satisfaction. On this basis, a distributed anytime-feasible resource allocation (DanyRA) algorithm is proposed. It is shown that the DanyRA algorithm converges to the exact optimal solution of DRAPs while preserving feasibility of the coupled inequality constraint from a feasible initialization. {\color{blue}To address constraint violations caused by external disturbances, we develop a recover-after-violation mechanism using a virtual queue with a minimum buffer, enabling feasibility recovery within a prescribed number of iterations.} We characterize the trade-off between convergence accuracy and violation robustness for a fixed positive buffer and establish exact convergence with a suitably decaying buffer. {\color{blue}DanyRA is further specialized to coupled equality constraints, with constraint preservation and an explicit linear convergence rate established.} Finally, a numerical example is provided for verification.
		\end{abstract}
		
	\end{frontmatter}

	\section{Introduction}
	\label{sec:introduction}
	Distributed resource allocation problem (DRAP) has received much attention due to its wide applications in various domains, including industrial
	Internet of Things \cite{2024_TII_Distributed_Multidomain_Resource_Allocation}, smart grids \cite{2025_Automatica_Privacy_preserving_distributed_optimization_for_economic_dispatch_in_smart_grids}, and multi-robot systems \cite{2021_COMST_Resource_Allocation_and_Service_Provisioning_in_Multi_Agent_Cloud_Robotics}. DRAP aims to optimally allocate limited resources among agents, meeting demand requirements under globally coupled constraints \cite{network_optimization_example}. These constraints are subject to physical limitations that may strictly prohibit violations. For example, in real-time control scenarios featuring fast-dynamics systems governed by inexpensive microcontrollers, inadequate processing capability may induce hard constraint violations, potentially leading to system failure \cite{2022_SCL_ROTEC,2024_TAC_Robust_to_Early_Termination_Model_Predictive_Control}. To mitigate such risks, ensuring constraint feasibility at all time steps during operation, i.e., anytime feasibility, is indispensable.
	
	Existing distributed algorithms from either the primal or the dual perspectives \cite{2023_Automatica_Distributed_delayed_dual_averaging_for_distributed_soptimization,2020_Automatica_Linear_convergence_of_primal_dual_gradient_methods,2019_arcontrol_A_survey_of_distributed_optimization} primarily focus on cost optimality,  typically satisfying constraints only asymptotically. These approaches are fundamentally insufficient for real-time systems, where deadline fulfillment and hard constraint inviolability are non-negotiable prerequisites \cite{Hard_Real_Time_Computing_Systems}. Practical deployments are further complicated by factors such as adversarial attacks or noise interference, which can induce unexpected disruptions \cite{2024_Automatica_Resilient_and_constrained_consensus_against_adversarial_attacks,2025_Automatica_Almost_sure_convergence_of_distributed_optimization_with_imperfect_information_sharing}. This requires that resource allocation mechanisms must not only ensure anytime feasibility under normal conditions but also exhibit robustness to fast restore constraint feasibility before deadlines when violations arise from interferences.
	The above discussions motivate the current study on designing distributed algorithms for DRAPs that incorporate these features.

	\subsection{Related works}
	For solving DRAPs with anytime feasibility, the pioneering work \cite{xiao2006optimal} proposes a weighted gradient descent method, where the weight of inter-agent gradient exchanges is designed to enforce anytime feasibility of coupled constraints. However, this work only considers scalar cost functions and the coupling within constraints is ignored. These  limitations are addressed by the algorithm in \cite{2021_CDC_A_New_Family_of_Feasible_Methods_for_Distributed_Resource_Allocation} based on barrier functions and the right-hand side allocation strategy  \cite{Nonlinear_Programming}. Nevertheless, this algorithm necessitates sharing private local cost functions. An extension that requires only gradient exchanges can be found in \cite{2023_Automatica_Distributed_safe_resource_allocation}. However, such approaches pose potential privacy-leakage risks, as gradients are sensitive information and their exchange can easily raise privacy concerns \cite{2017_NIPs_Collaborative_Deep_Learning_in_Fixed_Topology_Networks}. To tackle such issues, Ref. \cite{2019_TCNS_A_Dual_Splitting_Approach_for_Distributed_Resource_Allocation} designs a distributed algorithm from the dual perspective, ensuring feasibility solely via dual variable exchanges. Notably, these algorithms \cite{xiao2006optimal,2021_CDC_A_New_Family_of_Feasible_Methods_for_Distributed_Resource_Allocation,2023_Automatica_Distributed_safe_resource_allocation,2019_TCNS_A_Dual_Splitting_Approach_for_Distributed_Resource_Allocation} ensure feasibility by strictly enforcing the invariance of the (weighted) sum of decision variables across iterations. This surely requires proper initialization.
	Once they are improperly initialized or disruptions induce constraint violations,  they cannot regain constraint feasibility anymore.
	
	To address the above challenge, one algorithm involves  introducing safety margin by proactively tightening the original constraint to avoid  violation \cite{2025_TCNS_A_Safe_First_Order_Method_for_Pricing_Based_Resource_Allocation_in_Safety_Critical_Networks}. However, this strategy may fail when interferences are sufficiently large. That is, constraint violation robustness cannot be guaranteed. Another class of algorithms is developed based on primal decomposition, which reformulates DRAPs into decoupled subproblems \cite{Nonlinear_Programming}. By alternately applying dual averaging and a primal-dual solver to each subproblem, the algorithms proposed in \cite{2024_arxiv_Achieving_violation_free_distributed_optimization_under_coupling_constraints,2025_TAC_A_Continuous_Time_Violation_Free_Multi_Agent_Optimization_Algorithm}  ensure the anytime feasibility even under large interferences. 
	However, their requirement for exact solutions to convex subproblems at each iteration may incur prohibitive computational complexity. To mitigate this burden, control barrier functions (CBFs), a computationally efficient tool for the synthesis of safe controllers \cite{2017_TAC_Control_Barrier_Function_Based_Quadratic_Programs}, emerge as a potential solution. Leveraging CBFs and the connection between algorithm design for constrained problems and feedback controller design for nonlinear systems, a continuous-time safe gradient flow is proposed in \cite{2024_TAC_CBF_Based_Design_of_Gradient_Flows}, which achieves exponential decay of constraint violations
	but is not amenable to distributed implementation. It is further extended to a distributed algorithm in \cite{2023_CDC_Distributed_and_Anytime_Algorithm_for_Network_Optimization_Problems_with_Separable_Structure} by integrating  safe gradient flow with projected saddle-point dynamics. However, its regularization step induces approximate optimality, and notably, {\color{blue}its constraint violation decays to zero exponentially as time goes to infinity.} This results in constraint satisfaction missing the deadline, which remains unacceptable for real-time systems with deterministic deadlines and hard constraints. A critical open challenge remains: how to simultaneously achieve anytime feasibility guarantees, and constraint violation elimination before deadlines with low computational complexity.
	
	\subsection{Contribution}
	
	Our contributions are summarized as follows.
	
	1) Adopting the connection between optimization algorithms and feedback controllers designs, we develop a novel control-theoretic distributed anytime-feasible (DanyRA) algorithm  for solving DRAPs. Specifically,  based on primal decomposition techniques, DRAP is first decomposed into decoupled subproblems, with variables updated via a primal-dual method. We then adopt a control perspective for optimization algorithm design by treating the above derived variables as reference signals, the feasible set as a safe set, decision variables as states, and by synthesizing a safe feedback controller using CBFs. Moreover, theoretical analysis proves that the decision variables converge to the exact optimal solution while maintaining constraint feasibility at all time steps.
	
	2) {\color{blue}We develop a recover-after-violation mechanism by incorporating a virtual queue with a minimum buffer into DanyRA to provide violation robustness against interferences.} This enables immediate elimination of violations under small interferences. For sufficiently large ones, violations can be eliminated within an adjustable number of iterations, thereby ensuring the recovery and maintenance of feasibility before the deadlines, i.e., deadline-aware feasibility. Theoretical analysis shows its convergence to an arbitrarily small neighborhood of the optimal solution.  Notably, each agent only solves  a small-size quadratic program, offering generally better computational efficiency than \cite{2024_arxiv_Achieving_violation_free_distributed_optimization_under_coupling_constraints,2025_TAC_A_Continuous_Time_Violation_Free_Multi_Agent_Optimization_Algorithm}. 
	
	3) We characterize the fundamental trade-off between convergence accuracy and violation robustness for maintaining or recovering feasibility. To resolve this dilemma, we adopt a  virtual queue with decaying buffer that first achieves simultaneous  violation robustness and optimality to a certain extent.
	
	{\color{blue}4) For DRAPs with a coupled equality constraint, we further characterize the corresponding specialization of DanyRA and establish its equality-feasibility preservation property and linear convergence rate.}
	
	The rest of the paper is organized as follows. Section \ref{sec:problem_formulation} formulates the DRAP. The proposed algorithms are developed in Section \ref{sec:algorithm_development}. Theoretical analyses on the convergence and feasibility are given in Section \ref{sec:convergence_analysis}. The extension to DRAPs with a global equality constraint is provided in Section \ref{sec:extension}. The algorithms are numerically tested in Section \ref{sec:numerical_experiments} and Section \ref{sec:sonclusion} concludes the paper.
	
	\emph{Notations:} Vectors default to columns if not otherwise specified. Bold letter $\mathbf{x}\!\in\! \mathbb{R}^{n p}$ is defined as $\mathbf{x}\!=\!{[x_1^\intercal, \cdots, x_n^\intercal]}^\intercal$ with $x_i \!\in \!\mathbb{R}^{p}$. Gradient of a function $f$ at $x$ is $\nabla f(x)$. $\otimes$ denotes the Kronecker product. $\mathbf{1}$  is the vector with all one entries. $\hat{\mathbf{1}}=\mathbf{1}\otimes I$, where $I$ is the identity matrix with proper dimensions. For vectors, $\Vert \cdot \Vert$ denotes the 2-norm, \(|\cdot|\) is the element-wise absolute value operator.  $\succcurlyeq$  denotes the element-wise 	greater than or equal to. Similar notations are used for $\preccurlyeq$ and $\succ$. For matrices,  $\bar{\sigma}_X$ ($\underline{\sigma}_X$) denotes the maximum (minimum non-zero) eigenvalue of $X$ and $\kappa_X=\frac{\bar{\sigma}_X}{\underline{\sigma}_X}$. We use $blkdiag(X_1,\cdots,X_n)$ to denote the block-diagonal matrix with $X_1,\cdots,X_n$ as blocks.
	
	\section{Problem Formulation}
	\label{sec:problem_formulation}
	Consider the constraint-coupled DRAP with $n$ agents:
	\begin{align}
		\begin{split}
			& \min_{\mathbf{x}} \quad  f(\mathbf{x})=\sum\limits_{i=1}^n \ f_i(x_i)  \\
			&  \ \ \text{s.t.} \quad \sum_{i=1}^n A_ix_i\preccurlyeq\sum_{i=1}^nd_i,
		\end{split} \label{IDRAP}
	\end{align}
	where $f_i:\mathbb{R}^{p}\rightarrow\mathbb{R}$ is the agent $i$'s local convex cost function, $x_i\in \mathbb{R}^p$ is the local decision variable, {\color{blue}$d_i$ denotes the local contribution to the global resource budget $D$}, and $A_i \in \mathbb{R}^{m\times p} \ (p\geq m)$ is the coupling matrix with full row rank. {\color{blue}The decomposition of $D$ into \(\{d_i\}_{i=1}^{n}\) is introduced to facilitate the distributed reformulation below and does not change the global feasible set; any decomposition satisfying $\sum_{i=1}^{n}d_i=D$ can be adopted.} Throughout the paper, we assume problem (\ref{IDRAP}) has an optimal solution $\mathbf{x}^\star$.
	
	The communication network over which agents exchange information can be represented by an undirected graph $\mathcal{G}=\left( \mathcal{N},\mathcal{E} \right) $, where $\mathcal{N}\!=\!\{1,\cdots\!,n\}$ is the set of agents and $\mathcal{E}\!\subseteq\! \mathcal{N}\!\times\! \mathcal{N}$ denotes the set of edges, accompanied with a nonnegative weighted matrix $\mathcal{W}=\left[w_{i,j}\right]_{i,j\in\mathcal{N}}$. For any $ i,j \!\in\! \mathcal{N}$ in the network, $w_{i,j}\!>\!0$ denotes that agent $j$ can exchange information with agent $i$. The corresponding Laplacian matrix $\mathcal{L}=\left[l_{i,j}\right]_{i,j\in\mathcal{N}}$ is defined as $l_{i,i}=\sum_{j=1,j\neq i}^n{w_{i,j}}$, $l_{i,j}=-w_{i,j},\forall i,j\in\mathcal{N},i\neq j$.
	The collection of all individual agents that agent $i$ can communicate with is defined as its neighbors set $\mathcal{N}_i$.
	
	\begin{assumption}
		The undirected network $\ \mathcal{G}$ is connected and the weighted matrix $\mathcal{W}$ is doubly stochastic, i.e., $\mathcal{W}=\mathcal{W}^\intercal$,  $\mathcal{W}\mathbf{1}=\mathbf{1}$, and $\mathbf{1}^\intercal\mathcal{W}=\mathbf{1}^\intercal$.
		\label{connected}
	\end{assumption}
	Under Assumption \ref{connected}, we have $(\mathbf{1}^\intercal\mathcal{L})^\intercal\!=\!\mathcal{L}\mathbf{1}\!=\!\mathbf{0}$ and its second smallest eigenvalue is positive \cite{horn2012matrix}.
	
	Our main goal of the paper is to develop a distributed algorithm with low computational complexity for DRAP (\ref{IDRAP}), which ensures anytime feasibility under normal conditions and exhibits robustness to fast restore constraint feasibility upon violations.

	\section{Algorithm Development}
	\label{sec:algorithm_development}
	In this section, we develop a distributed anytime-feasible resource allocation algorithm for DRAP (\ref{IDRAP}).
	
	The constraint in DRAP (\ref{IDRAP}) is coupled, involving all local variables, and thus is typically inaccessible to individual agents in practice. Moreover, direct global projection onto the feasible set is inapplicable under such a coupled constraint. To address this, we first derive an equivalent formulation of DRAP (\ref{IDRAP}) to resolve the coupling among agents' decision variables.
	
	\emph{Proposition 1}: Under Assumption \ref{connected}, $\mathbf{x}^\star$ is an optimal solution of (\ref{IDRAP}) if and only if there are $\mathbf{y}^\star\in \mathbb{R}^{nm}$ and $\boldsymbol{\delta}^\star\geq0$ such that $(\mathbf{x}^\star,\mathbf{y}^\star,\boldsymbol{\delta}^\star)$ is an optimal solution of the following optimization problem
	\begin{align}
		\begin{split}
			& \min_{\boldsymbol{\delta}\succcurlyeq0,\mathbf{x},\mathbf{y}} \qquad  f(\mathbf{x})=\sum\limits_{i=1}^n \ f_i(x_i)  \\
			&  \ \text{s.t.} \ A_ix_i\!+\!\delta_i\!-\!d_i\!+\!\sum\nolimits_{j\in\mathcal{N}_i}\!w_{i,j}(y_i\!-\!y_j)\!=\!0, \forall i\!\in\! \mathcal{N},
		\end{split} \label{IDRAP_1}
	\end{align}
	{\color{blue}where \(y_i\) is an auxiliary variable and \(\delta_i\succeq0\) is the slack variable associated with the local constraint of agent \(i\).}
	
	\emph{Proof}: See Appendix \ref{Proof_of_Proposition_1}. $\hfill \blacksquare$
	
	{\color{blue}The role of \(y_i\) is to enable the coupled constraint in (1) to be represented by local constraints involving neighboring information. Under Assumption~1, the Laplacian terms satisfy $\sum_{i=1}^{n}\sum_{j\in\mathcal N_i}w_{i,j}(y_i-y_j)=0$, and hence do not alter the original global coupled constraint. Therefore, together with the variable \(\delta_i\), formulation \eqref{IDRAP_1} provides an equivalent distributed representation of problem \eqref{IDRAP}.}
	
	According to Karash-Kuhn-Tucker conditions \cite{1970_Convex_Analysis}, the optimal solution of DRAP (\ref{IDRAP_1}) is attained if and only if the following conditions are satisfied:
	\begin{subequations}
		\label{KKT_inequality}
		\begin{align}
			& \nabla f(\mathbf{x}^\star)+\mathbf{A}^\intercal\boldsymbol{\lambda}^\star=0, \ \boldsymbol{\mathcal{L}}\boldsymbol{\lambda}^\star=0,  \label{KKT_inequality_1} \\ 
			& \boldsymbol{\delta}^\star\succcurlyeq 0, \ {\boldsymbol{\lambda}^\star}^\intercal{\boldsymbol{\delta}^\star}=0, \ \mathbf{A}\mathbf{x}^\star+\boldsymbol{\mathcal{L}}\mathbf{y}^\star+\boldsymbol{\delta}^\star-\mathbf{d}=0, \label{KKT_inequality_2} 
		\end{align}
	\end{subequations}
	where $\mathbf{A}=blkdiag(A_1,\cdots,A_n)$, $\boldsymbol{\mathcal{L}}=\mathcal{L}\otimes I_m$, and $\boldsymbol{\lambda}$ is the Lagrange multiplier for the constraints of (\ref{IDRAP_1}).
	
	Motivated by the feedback-control design of Lagrange multipliers \cite{2025_TAC_A_New_Framework_for_Constrained_Optimization_via_Feedback_Control_of_Lagrange_Multipliers}, we construct a distributed primal--dual iteration for problem \eqref{IDRAP_1}. The multiplier update combines integral accumulation of the constraint residual with an additional feedback term based on the primal stationarity residual, giving
	\begin{subequations}
		\label{primal_dual}
		\begin{align}
			x^\prime_{i,k+1}&=x^\prime_{i,k}\!-\!\alpha(\nabla f_i(x^\prime_{{i,k}})\!+\!A_i^\intercal(z_{i,k}\!-\!d_i\!+\!\lambda_{i,k})),\\
			y_{i,k+1}
			&= y_{i,k}
			-\alpha \sum_{j\in\mathcal{N}_i} w_{i,j}
			\Bigl[
			(z_{i,k}-d_i+\lambda_{i,k}) \notag\\
			&\qquad\qquad
			-(z_{j,k}-d_j+\lambda_{j,k})
			\Bigr],\\
			\delta_{i,k+1}&=\max\{\delta_{i,k}-\alpha({z}_{i,k}-d_i+{\lambda}_{i,k}),0\},\\	z_{i,k+1}&=\!A_ix^\prime_{i,k\!+\!1}\!+\!\!\sum_{j\in\mathcal{N}_i}\!\!w_{i,j}(y_{i,k\!+\!1}\!-\!y_{j,k\!+\!1})\!+\!\delta_{i,k\!+\!1}, \\
			\lambda_{i,k+1}&\!=\!\lambda_{i,k}\!+\!\beta\left(z_{i,k+1}\!-\!d_i\!-\!\eta A_i(A_i^\intercal\lambda_{i,k}\!\!+\!\!\nabla f_i(x^\prime_{{i,k}}))\right),
		\end{align}
	\end{subequations}
	where $x'_i$ is the virtual decision variable, $\lambda_i$ is the dual variable, and $z_i$ is the intermediate variable. {\color{blue}In the multiplier update, $z_{i,k+1}-d_i$ provides the constraint-residual feedback, while the additional term
	$-\eta A_i(A_i^\intercal\lambda_{i,k}+\nabla f_i(x'_{i,k}))$
	incorporates the primal stationarity residual and vanishes at a KKT point.} In this way, $\mathbf{x}^\prime$, $\mathbf{y}$, and $\boldsymbol{\delta}$ converge asymptotically to an optimal solution $(\mathbf{x}^\star,\mathbf{y}^\star,\boldsymbol{\delta}^\star)$ of problem (\ref{IDRAP_1}) \cite{2025_TAC_A_New_Framework_for_Constrained_Optimization_via_Feedback_Control_of_Lagrange_Multipliers}.
	
	To further achieve anytime feasibility, we treat the variable $\mathbf{x}^\prime$ derived by (\ref{primal_dual}) as the nominal state and construct the following optimization problem 
	\begin{align}
		\begin{split}
			& \min_{x_i} \quad  g_i(x_{i})=\frac{1}{2}\Vert x_{i}-x^\prime_{i,k+1}\Vert^2  \\
			&  \ \ \text{s.t.} \quad A_ix_i\!+\!\delta_i^\star\!-\!d_i\!+\!\!\sum\limits_{j\in\mathcal{N}_i}\!\!w_{i,j}(y_i^\star\!-\!y_j^\star)\!=\!0
		\end{split} \label{problem_before_CBF}
	\end{align}
	for agent $i$ to track the nominal state at each time $k$, while referring to the constraint.
	We can adopt the penalty method to solve the problem (\ref{problem_before_CBF}), resulting in the following iterative process
	\begin{align}
		x_{i,k+1}=x_{i,k}-\nabla g_i(x_{i,k})-\nabla h_i(x_{i,k})^\intercal u_{i,k}, \label{control_system}
	\end{align}
	where $h_i(x_{i})\!=\!A_ix_i\!+\!\delta_i^\star\!-\!d_i\!+\!\!\sum\limits_{j\in\mathcal{N}_i}\!\!\!w_{i,j}(y_i^\star\!-\!y_j^\star)$, and $u_{i,k}$ is the  penalty parameter that should be properly designed \cite{Nonlinear_Programming}. 
	
	Under Assumption \ref{connected}, it can be verified that $x_{i,k+1}$ satisfies the coupled constraint in the DRAP(\ref{IDRAP}), i.e, it is a feasible solution, when $h_i(x_{i,k+1})=0$ for all $i\in\mathcal{N}$ and $k\geq0$. This motivates us to define 
	\begin{align}
		\mathcal{C}_i=\{x_i\in\mathbb{R}^p|h_i(x_{i})=0\} \label{safe_set}
	\end{align} 
	as a safe set. Since $\nabla h_i(x_{i})\nabla h_i(x_{i})^\intercal=A_iA_i^\intercal$ is full-rank given that $A_i$ has full row rank, the following set 
	\begin{align}
		\mathcal{K}_i(x_i)=&\big\{u_i\in\mathbb{R}^m| \notag \\
		&-\nabla h_i\nabla h_i^\intercal u_i=\nabla h_i \nabla g_i(x_i) -\gamma h_i(x_{i})\big\} \label{admissible_set}
	\end{align}
	is nonempty  $\forall x_i\in\mathbb{R}^p$, where the constant $0<\gamma<1$.
	
	According to \cite{2024_TAC_CBF_Based_Design_of_Gradient_Flows,2017_TAC_Control_Barrier_Function_Based_Quadratic_Programs}, relations (\ref{safe_set}) and (\ref{admissible_set}) imply that $h_i$ is a CBF of $\mathcal{C}_i$ for the control system (\ref{control_system}) with $u_{i,k}$ regarded as the control input. Based on this observation, we can now utilize the CBF technique to design the penalty parameter $u_i$ by synthesizing the safe controller as follows
	\begin{align}
		u_{i,k}&\in {\arg\min}_{u_i\in\mathcal{K}_i(x_{i,k})}\left\{\Vert\nabla h_i(x_{i,k})^\intercal u_i\Vert^2\right\}. \label{safe_feedback_controller}
	\end{align}
	{\color{blue}Since the function $h_i(x_i)$ is affine, substituting  \eqref{safe_feedback_controller} into \eqref{control_system} yields
		$h_i(x_{i,k+1})=(1-\gamma)h_i(x_{i,k})$.
		Thus, the given update preserves $\mathcal{C}_i$, i.e., once $x_{i,k}\in\mathcal{C}_i$, $x_{i,k+t}\in\mathcal{C}_i$ holds $\forall t\geq0$, and contracts nonzero residuals geometrically.} 
	
\begin{figure}[thpb]
		\centering
		\includegraphics[height=3.5cm]{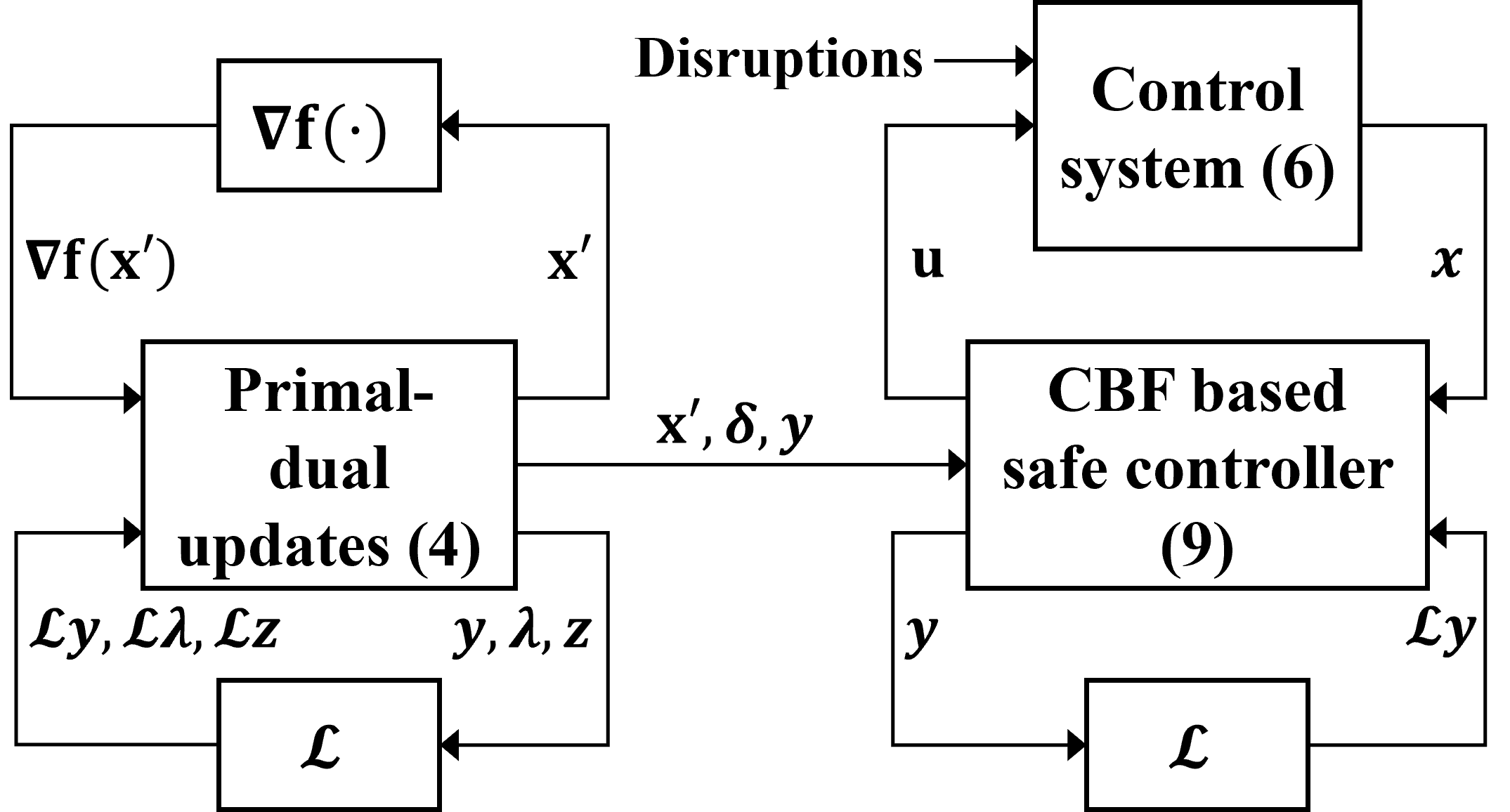}
		\caption{Schematic view of the overall safe  controller design.}
		\label{algorithm_design}
	\end{figure}

	One main issue is, however, that $\delta_i^\star$ and $y_i^\star$ within $h_i$ are unknown.  To have an amenable implementation, in the iterative process, $\delta_{i,k+1}$ and $y_{i,k+1}$ are employed instead. Combining (\ref{control_system}) with (\ref{admissible_set}) and (\ref{safe_feedback_controller}), the update rule of $x_i$ can be written as
	\begin{align}
		x_{i,k+1}&\in {\arg\min}_{x_i\in\mathcal{K}^x_i(x_{i,k})}\left\{\Vert x_i-x^\prime_{i,k+1}\Vert^2\right\}, \label{safe_x_update}
	\end{align}
	where $
	\mathcal{K}^{x}_i(x_{i,k})=\{x_i\in\mathbb{R}^p| A_i x_i=A_i x_{i,k} -\gamma \big(A_ix_{i,k}\!+\!\delta_{i,k+1}\!-\!d_i\!+\!\sum\nolimits_{j\in\mathcal{N}_i}w_{i,j}(y_{i,k+1}\!-\!y_{j,k+1})\big)\!+\!(1\!-\!\gamma)(\delta_{i,k}\!-\delta_{i,k+1})\}$. It is noted that the set $\mathcal{K}^x_i(x_{i,k})$ is nonempty since $A_i$ has full row rank. {\color{blue} This construction adapts CBFs to discrete-time distributed optimization by using the equivalent formulation \eqref{IDRAP_1} to enforce the coupled constraint through local updates and neighboring information.}
	Moreover, in set $\mathcal{K}^{x}_i(x_{i,k})$ we introduce the term $(1-\gamma)(\delta_{i,k}\!-\delta_{i,k+1})$ to compensate the constraint violation induced by replacing $\delta_i^\star$ with $\delta_{i,k+1}$.  {\color{blue}In contrast, no compensation is required for \(y_i\), since its Laplacian contribution satisfies \(\sum_{i=1}^n\sum_{j\in\mathcal N_i}w_{i,j}(y_i-y_j)=0\) under Assumption \ref{connected} and thus does not affect the global coupled constraint.} {\color{blue}Moreover, \eqref{safe_x_update} minimally corrects the nominal primal--dual iterate, while subsequent primal--dual updates in \eqref{primal_dual} continue to use $\mathbf{x}'_k$, preserving the nominal dynamics.  Moreover, the induced discrepancy is explicitly accounted for in the convergence analysis in Section~\ref{sec:convergence_analysis}.} See Fig. \ref{algorithm_design} for a schematic view of the overall safe controller design.
	
	A distributed anytime-feasible resource allocation (DanyRA) algorithm is  summarized in Algorithm \ref{alg2}, 
	

	\begin{algorithm}
		\renewcommand{\algorithmicrequire}{\textbf{Input:}}
		\footnotesize
		\caption{\bf : DanyRA algorithm}
		\label{alg2}
		\begin{algorithmic}[1]
			\STATE {Parameters: 
				proper stepsizes~$\alpha,\beta,\eta,\gamma>0$; $\omega\geq0$\\
				Initialization: Arbitrary $x^\prime_{i,0},x_{i,0}\in \mathbb{R}^{p},\lambda_{i,0}\in \mathbb{R}^{m}, y_{i,0}=\delta_{i,0}=0$.}
			\FOR{$k=0, 1, \cdots, K-1 $}
			\STATE \emph{Information exchange between agents:}  \\
			$\bar{\lambda}_{i,k}\!=\!\sum\limits_{j\in\mathcal{N}_i}w_{i,j}(\lambda_{i,k}\!-\!\lambda_{j,k}),$ $\bar{y}_{i,k}\!=\!\sum\limits_{j\in\mathcal{N}_i}\!w_{i,j}(y_{i,k}\!-\!y_{j,k}),$\\ 
			$\bar{z}_{i,k}=\sum\limits_{j\in\mathcal{N}_i}w_{i,j}((z_{i,k}-d_i)-(z_{j,k}-d_j))$,\\
			\STATE \emph{Intermediate variable update:}  \\
			$x^\prime_{i,k+1}=x^\prime_{i,k}-\alpha(\nabla f_i(x^\prime_{{i,k}})+A_i^\intercal(z_{i,k}-d_i+\lambda_{i,k}))$,\\
			\STATE \emph{Virtual queue update with minimum buffer:} \\
			$ y_{i,k+1}=y_{i,k}-\alpha(\bar{z}_{i,k}+\bar{\lambda}_{i,k})$,\\
			$ \delta_{i,k+1}=\max\{\delta_{i,k}-\alpha({z}_{i,k}-d_i+{\lambda}_{i,k}),{\omega}\}$,\\		
			\STATE \emph{Dual variable update:} \\
			$\bar{y}_{i,k+1}=\sum\limits_{j\in\mathcal{N}_i}w_{i,j}(y_{i,k+1}-y_{j,k+1}),$\\
			$z_{i,k+1}=A_ix^\prime_{i,k+1}+\bar{y}_{i,k+1}+\delta_{i,k+1}$, \\
			$\lambda_{i,k+1}=\lambda_{i,k}+\beta\big(z_{i,k+1}-d_i-\eta A_i(A_i^\intercal\lambda_{i,k}+\nabla f_i(x^\prime_{{i,k}}))\big)$,\\	
			\STATE \emph{Decision variable update:}  \\
			$ x_{i,k+1}=\arg\min_{x_i\in\mathcal{K}^{x}_i(x_{i,k})} \frac{1}{2}\Vert x_i-x^\prime_{i,k+1}\Vert^2.$
			\ENDFOR
		\end{algorithmic}
	\end{algorithm}
	
	For the DanyRA algorithm, to tackle potential constraint violations induced by unexpected disruptions, a virtual queue with minimum buffer $\omega$ is further introduced, see step $5$. The update is related to dual variable update and virtual queue update \cite{Nonlinear_Programming,Stochastic_Network_Optimization_with_Application_to_Communication_and_Queueing_Systems}. An important difference is that a minimum buffer $\omega>0$ is added rather than $\omega=0$. This buffer serves as a cushion: once exceeded, a restoring force is applied to pull the decision back, thereby enabling rapid restoration of constraint satisfaction. 
	This addresses the limitation of existing methods \cite{2024_TAC_CBF_Based_Design_of_Gradient_Flows,2019_arcontrol_A_survey_of_distributed_optimization}, which only achieve asymptotic recovery of violated global constraints. 
	
	{\color{blue}Traditional robust CBF methods typically incorporate a known disturbance bound to preserve safety throughout the evolution of physical systems. Here, the states in the control interpretation are optimization iterates and do not directly represent physical quantities. For settings permitting temporary violations, our recover-after-violation mechanism responds to the realized violation and restores feasibility through the buffer $\omega$ and parameter $\gamma$, without requiring an a priori disturbance bound. This avoids reliance on potentially conservative worst-case bounds and is particularly useful when disturbance bounds are unknown or overestimated. The finite-time recovery guarantee is established in Section~\ref{sec:convergence_analysis}. As a final remark, each agent corrects its primal–dual iterate by solving a local quadratic program, even when its local cost function is nonquadratic. This simplifies the local computation compared with the general convex subproblems \cite{2024_arxiv_Achieving_violation_free_distributed_optimization_under_coupling_constraints,2025_TAC_A_Continuous_Time_Violation_Free_Multi_Agent_Optimization_Algorithm}, as numerically verified in Section \ref{sec:numerical_experiments}.}
	
	\section{Convergence Analysis}
	\label{sec:convergence_analysis}
	In this section, we first prove the convergence properties including optimality and anytime feasibility of the DanyRA algorithm. Then, we characterize the trade-off between convergence accuracy and violation robustness 
	
	We make the following assumptions on DRAP (\ref{IDRAP}):
	
	\begin{assumption}
		Each cost function $f_i$ is $\ell$-Lipschitz smooth, i.e., there is some $\ell < \infty$ such that
		$\Vert \nabla f_i(x)-\nabla f_i(x^\prime)\Vert \leq \ell\Vert x-x^\prime\Vert, \ \forall x, x^\prime \in \mathbb{R}^p$.
		\label{Lipschitz_smooth}
	\end{assumption}
	
	\begin{assumption}
		Each cost function $f_i$ is restricted $\mu$-strongly convex with respect to the optimal solution $x^\star$, i.e., there is some $\mu >0$ such that $
		(x-x^\star)^\intercal \left(\nabla f_i(x)-\nabla f_i(x^\star)\right)\geq \mu \Vert x-x^\star\Vert^2, \ \forall x\in \mathbb{R}^p$.
		\label{strong_convexity}
	\end{assumption}
	Assumption \ref{Lipschitz_smooth} is standard in related works \cite{2019_arcontrol_A_survey_of_distributed_optimization}. The restricted $\mu$-strongly convexity \cite{2015_siamjopt_Extra} in Assumption \ref{strong_convexity} is weaker than the strongly convexity that is
	commonly used to derive the linear convergence rate of gradient based methods.
	
	\subsection{Convergence Analysis for DanyRA Algorithm}
	We first establish relations between the fixed point of DanyRA algorithm and the optimal solution of DRAP (\ref{IDRAP_1}). As given below, they are consistent when $\omega=0$.
	
	\emph{Proposition 2}: Under Assumption \ref{connected}, any fixed point of DanyRA algorithm with $\omega=0$, denoted as $(\hat{\mathbf{x}},\hat{\mathbf{y}},\hat{\boldsymbol{\delta}},\hat{\boldsymbol{\lambda}})$ satisfies the optimality conditions
	given in (\ref{KKT_inequality}).
	
	\emph{Proof}: See Appendix \ref{Proof_of_Proposition_2}.  
	$\hfill \blacksquare$
	
	{\color{blue}The proof also shows that the virtual decision is admissible for the safety-correction subproblem at the fixed point. The minimum-distance objective therefore gives $\hat{\mathbf{x}}=\hat{\mathbf{x}}'$, so the correction introduces no discrepancy between the actual and virtual decisions at the KKT point.} For a positive $\omega$, we obtain the following result.
	
	\emph{Proposition 3}: Under Assumption \ref{connected}, any fixed point of the DanyRA algorithm with $\omega>0$, denoted as $(\hat{\mathbf{x}},\hat{\mathbf{y}},\hat{\boldsymbol{\delta}},\hat{\boldsymbol{\lambda}})$ is a feasible solution of (\ref{IDRAP}) and it satisfies
	\begin{align}
		\Vert\hat{\mathbf{x}}-\mathbf{x}^\star\Vert \leq \frac{\ell\sqrt{n}}{\mu\underline{\sigma}_\mathbf{A}} \omega. \label{Sensitivity_result}
	\end{align}
	\emph{Proof}: See Appendix \ref{Proof_of_Proposition_3}.  
	$\hfill \blacksquare$
	
	For convergence analysis, we first define the virtual decision error term $\tilde{\mathbf{x}}^\prime_{k}=\mathbf{x}^\prime_{k}-\hat{\mathbf{x}}^\prime$, dual error term $\tilde{\boldsymbol{\lambda}}_{k}=\boldsymbol{\lambda}_k-\hat{\boldsymbol{\lambda}}$, and auxiliary error terms $\tilde{\mathbf{y}}_{k}=\mathbf{y}_{k}-\hat{\mathbf{y}}$,  $\tilde{\boldsymbol{\delta}}_{k}=\boldsymbol{\delta}_k-\hat{\boldsymbol{\delta}}$. 
	
	Moreover, we introduce the metric
	\begin{align}
		V_{k}=&\Vert \tilde{\mathbf{x}}^\prime_{k} \Vert^2+\Vert\tilde{\mathbf{y}}_{k} \Vert^2+\Vert\tilde{\boldsymbol{\delta}}_{k}\Vert^2+\frac{\alpha}{\beta}\Vert \tilde{\boldsymbol{\lambda}}_{k} \Vert^2 \notag \\
		&+\frac{\alpha(1\!-\!3\beta)}{2}\Vert\tilde{\mathbf{z}}_k\Vert^2\!+\!\frac{\underline{\sigma}_\mathbf{A}^2\alpha(1\!-\!3\beta)}{8\gamma^2}\Vert\mathbf{x}_{k}\!-\!\mathbf{x}^\prime_{k}\Vert^2,
	\end{align}
	where  $\mathbf{z}_{k}=\mathbf{A}\mathbf{x}^\prime_k+\boldsymbol{\mathcal{L}}\mathbf{y}_k+\boldsymbol{\delta}_k$, $\hat{\mathbf{z}}=\mathbf{A}\hat{\mathbf{x}}+\boldsymbol{\mathcal{L}}\hat{\mathbf{y}}+\hat{\boldsymbol{\delta}}$ and $\tilde{\mathbf{z}}_{k}=\mathbf{z}_{k}-\hat{\mathbf{z}}$. The first four terms of $V_{k}$ correspond to the virtual decision, auxiliary and dual errors, and the last term quantifies the decision tracking error. {\color{blue}The fifth term measures the residual of the	constraint in \eqref{IDRAP_1} evaluated at the virtual	decision $\mathbf{x}'_k$. Since the corresponding residual at the actual decision $\mathbf{x}_k$ equals $\tilde{\mathbf{z}}_k+\mathbf{A}(\mathbf{x}_k-\mathbf{x}'_k)$, the fifth and last terms jointly control this residual.}
	
	\emph{Lemma 1:}  Under Assumptions \ref{connected}-\ref{strong_convexity}, the variables generated by DanyRA algorithm satisfy
	\begin{align}
		&V_{k+1}-V_k\leq \frac{6\alpha\bar{\sigma}_\mathbf{A}^2(1\!+\!3c)\!-\!1}{2}\Vert\Delta\mathbf{x}^\prime_{k\!+\!1}\Vert^2\! \notag\\
		&+\!(3\alpha\bar{\sigma}_\mathcal{L}^2(1\!+\!4c)\!-\!1)\Vert\Delta\mathbf{y}_{k\!+\!1}\Vert^2\!\!+\!(3\alpha(1\!+\!3c)\!-\!1)\Vert\Delta\boldsymbol{\delta}_{k\!+\!1}\Vert^2  \notag \\
		&+\! \alpha(2\ell^2\alpha\!-\!2\mu\!+\!\eta\ell^2(3\beta\eta\!+\!1))\Vert\tilde{\mathbf{x}}^\prime_k\Vert^2\notag \\
		&+\alpha\eta(3\beta\eta\bar{\sigma}_\mathbf{A}^2-\underline{\sigma}_\mathbf{A}^2)\Vert\tilde{\boldsymbol{\lambda}}_k\Vert^2+\frac{\alpha(3\beta-1)}{2}\Vert\tilde{\mathbf{z}}_k\Vert^2\notag \\
		&+\frac{\underline{\sigma}_\mathbf{A}^2\alpha(1\!-\!3\beta)}{8\gamma^2}\left({4(1-\gamma)^2\kappa_\mathbf{A}^2}\!-\!1\right)\Vert\mathbf{x}_k-\mathbf{x}^\prime_k\Vert^2, \label{result_Lemma_1}
	\end{align}
	where $\Delta\mathbf{x}_{k+1}=\mathbf{x}_{k+1}-\mathbf{x}_{k}$, $\Delta\mathbf{y}_{k+1}=\mathbf{y}_{k+1}-\mathbf{y}_{k}$, $\Delta\boldsymbol{\delta}_{k+1}=\boldsymbol{\delta}_{k+1}-\boldsymbol{\delta}_{k}$, and the constant $c=\frac{(1-\gamma)^2(1-3\beta)}{2\gamma^2}$.
	
	\emph{Proof:} See Appendix \ref{Proof_of_Lemma_1}.
	$\hfill \blacksquare$
	
	Then, the convergence properties of DanyRA algorithm are given as follows.
	
	\emph{Theorem 1:} Under Assumptions \ref{connected}-\ref{strong_convexity}, if the following conditions  are satisfied
	\begin{subequations}
		\label{conditions_Theorem_1}
		\begin{align}		&\alpha<\min\left\{\frac{1}{6\bar{\sigma}_{\mathbf{A}}^2(1+3c)},\frac{1}{3\bar{\sigma}_{\mathcal{L}}^2(1+4c)},\frac{1}{6(1+3c)}, \right. \notag \\
			&\quad \ \ \left.\frac{2\mu\!-\!\eta\ell^2(3\beta\eta\!+\!1)}{2\ell^2},2\eta(\underline{\sigma}_\mathbf{A}^2\!-\!3\beta\eta\bar{\sigma}_\mathbf{A}^2),1\!-\!3\beta\right\}, \\
			& \beta\!<\!\min\!\left\{\frac{1}{3},\frac{1}{3\eta}\!\Big(\frac{2\mu}{\eta\ell^2}\!-\!1\Big), \frac{1}{3\eta\kappa_\mathbf{A}^2}\right\}, \\ &\eta\!<\!\frac{2\mu}{\ell^2},  1\!-\!\frac{1}{2\kappa_\mathbf{A}}\!<\!\gamma\!<\!1,
		\end{align}
	\end{subequations}
	then, the proposed DanyRA algorithm satisfies
	
	\romannumeral1) $\{\mathbf{x}_{k}\}_{k\geq0}$ converges to the neighborhood of the optimal solution of (\ref{IDRAP}), and the convergence error is bounded by
	\begin{align}
		\lim_{k\rightarrow\infty}\Vert \mathbf{x}_{k}-\mathbf{x}^\star \Vert^2\leq\frac{\ell\sqrt{n}}{\mu\underline{\sigma}_\mathbf{A}} \omega; \label{Theorem_1_result_1}
	\end{align}
	\romannumeral2) when the global constraint is violated for some $k_0\geq 1$, the violation $C^{vio}_{k_0}=\Vert\max\{\sum_{i=1}^n (A_ix_{i,k_0}-d_{i}),0\}\Vert_\infty$ decreases to zero in a finite number of iterations. Specifically, the violation remains zero after $t$ iterations, where
	\begin{align}
		t\geq \frac{\ln\Big(\frac{n\omega}{C^{vio}_{k_0}+C_\delta}\Big)}{\ln(1-\gamma)}, \label{Theorem_1_result_4}
	\end{align}
	{  where the constant $C_\delta=\sqrt n(\|\hat{\boldsymbol\delta}\|+\sqrt{V_0})$;}
	
	\romannumeral3) once the relation $\sum_{i=1}^n \left(A_ix_{i,k_0}+\delta_{i,k_0}-d_{i}\right)\preccurlyeq \frac{n\omega\mathbf{1}}{1-\gamma}$ holds for some $k_0\geq 1$, the global  constraint remains satisfied thereafter:
	\begin{align}
		\sum_{i=1}^n A_ix_{i,k_0+t}\preccurlyeq\sum_{i=1}^nd_{i}, \ \forall t\geq 1.
	\end{align}
	\emph{Proof:} See Appendix \ref{Proof_of_Theorem_1}.
	$\hfill \blacksquare$
	
	{\color{blue}The admissible parameter set in \eqref{conditions_Theorem_1} is nonempty. First, choose $0<\eta<2\mu/\ell^2$ and select $\gamma$ from its prescribed interval, which is nonempty since $\kappa_{\mathbf A}$ is finite. The choice of $\eta$ makes all upper bounds on $\beta$ positive, so one can choose $\beta>0$ below these bounds. These choices make $c$  nonnegative and all upper bounds on $\alpha$ positive, so a sufficiently small $\alpha>0$ is admissible.}
	
	{\color{blue}\emph{Remark}: The recovery result in Theorem 1 assumes that no new disruption is introduced during recovery. It applies to intermittent disruptions when the inter-disruption intervals are sufficiently long for recovery to be completed. Persistent or high-frequency disturbances may introduce new violations before the ongoing recovery process is completed and therefore require additional assumptions on their magnitude and frequency. Extending the proposed framework to such disturbance settings is left for future work.}
	
	It is shown in \romannumeral1) that $\mathbf{x}_k$ converges to an arbitrarily small neighborhood of the optimal solution as $\omega\rightarrow0$.
	Moreover, from \romannumeral3), it follows that with a feasible initialization, the derived decisions remain feasible at every iteration, i.e., anytime feasibility is achieved. Statement in \romannumeral3) also shows that $\omega$ provides robustness against constraint violations that is less than $\frac{n\omega}{1-\gamma}$. Specifically, a small violation of the constraint in (\ref{IDRAP_1}) does not influence the feasibility of the derived solution thereafter for the original DRAP (\ref{IDRAP}). Moreover,  \romannumeral2) indicates that a large $\omega$ enhances the DanyRA algorithm’s ability to eliminate constraint violations rapidly with fewer iterations.  
	
	{\color{blue}Theorem 1 gives an explicit bound on the number of iterations required to restore feasibility, which guides the choice of \(\omega\) and \(\gamma\) for a prescribed recovery deadline. For a fixed $\omega$, increasing $\gamma$ within its admissible range reduces the contraction factor $1-\gamma$ and thereby accelerates feasibility recovery, without changing the accuracy bound in \eqref{Theorem_1_result_1}. Its influence on the optimality convergence rate, however, is coupled with $\alpha$, $\beta$, and $\eta$ through \eqref{conditions_Theorem_1}, and thus cannot be characterized by $\gamma$ alone.}

	\subsection{Trade-off between Convergence Accuracy and Violation Robustness}
	{\color{blue}A larger buffer $\omega$ enhances feasibility recovery but enlarges the steady-state error bound in \eqref{Theorem_1_result_1}. Thus, a fundamental trade-off emerges between convergence accuracy and violation robustness, which can be characterized by selecting the smallest admissible $\omega$ that guarantees recovery within a prescribed horizon $\underline t$ according to \eqref{Theorem_1_result_4}. Substituting this choice into \eqref{Theorem_1_result_1} yields the tightest corresponding accuracy bound:}
	\begin{align}
		\lim_{k\rightarrow\infty}\Vert \mathbf{x}_{k}-\mathbf{x}^\star \Vert^2\leq\frac{\ell(1-\gamma)^{\underline{t}+1}{(C^{vio}_{k_0}+C_\delta)}}{\mu\underline{\sigma}_\mathbf{A}\sqrt{n}}, \label{tradeoff}
	\end{align}
	where the integer $\underline{t}>0$ is the smallest $t$ satisfying (\ref{Theorem_1_result_4}). This indicates that when $\gamma$ is fixed, rapid  feasibility recovery (i.e., small $\underline{t}$) results in a large convergence error.
	
	As previously discussed, a small constraint violation can be immediately eliminated via one iteration, where the threshold $\frac{n\omega}{1-\gamma}$ is also dependent on $\gamma$. {\color{blue}Given a bounded violation magnitude denoted by $\bar C^{vio}$, we can select a sufficiently small $1-\gamma$ to satisfy $\bar C^{\rm vio}\leq \frac{n\omega}{1-\gamma}-C_\delta$, thereby ensuring feasibility recovery via one iteration under the conditions of Theorem 1.}
	
	From Theorem 1, the DanyRA algorithm with a virtual queue with minimum buffer achieves exact convergence only when $\omega = 0$, but this requires sacrificing robustness against constraint violations.  One potential solution is employing a virtual queue with decaying buffer $\omega_{k}$.
	
	The convergence properties of DanyRA algorithm with decaying parameter $\omega_k$ are given as follows.
	
	\emph{Theorem 2:} Under Assumptions \ref{connected}-\ref{strong_convexity}, if the conditions in Theorem 1 and $\sum\omega_{k}^2<\infty$ hold, the proposed DanyRA algorithm satisfies
	
	\romannumeral1) $\{\mathbf{x}_{k}\}_{k\geq0}$ converges to the exact optimal solution of DRAP (\ref{IDRAP}), i.e., $
	\lim_{k\rightarrow\infty}\Vert \mathbf{x}_{k}-\mathbf{x}^\star \Vert^2=0$;
	
	\romannumeral2) when the global constraint is violated for some $k_0\geq 1$, the violation $C^{vio}_{k_0}$ decreases to zero in a finite number of iterations, i.e., the violation remains zero for all $k\geq k_0+t$ iterations, where $t$ satisfies
	\begin{align}
		{(1-\gamma)^{t} (C^{vio}_{k_0}+C_\delta^\prime)\leq n \omega_{k_0+t}, \label{Theorem_3_result_3}}
	\end{align}
	{  where the constant $C_\delta^\prime=\sqrt n(\|\hat{\boldsymbol\delta}\|+\sqrt{V_0+6n\sum_{j=0}^\infty\omega_j^2})$;}
	
	\romannumeral3) once the relation $\sum_{i=1}^n \left(A_ix_{i,k_0}\!+\!\delta_{i,k_0}\!-\!d_{i}\right)\!\preccurlyeq\! \frac{n\omega_{k_0+1}\mathbf{1}}{1-\gamma}$ holds for some $k_0\geq 1$, the global  constraint remains satisfied thereafter:
	\begin{align}
		\sum_{i=1}^n A_ix_{i,k_0+t}\preccurlyeq\sum_{i=1}^nd_{i}, \ \forall t\geq 1.
	\end{align}
	\emph{Proof:} See Appendix \ref{Proof_of_Theorem_2}.
	$\hfill \blacksquare$
	
	With a properly chosen $\omega_{k}$, such as $\omega_{k}=1/(k+1),\forall k\geq0$, the condition (\ref{Theorem_3_result_3}) can always be feasible when $t\geq\underline{t}$ for some $\underline{t}\geq1$.
	In this way, exact convergence to the optimal solution of DRAP (\ref{IDRAP}), anytime feasibility, and robustness to constraint violation are simultaneously ensured. However, it should be noted that the robustness to constraint violations diminishes as $k$ increases. This may be acceptable since the decision variables essentially stabilize after a certain number of iterations. Users can design $\omega_k$ and the parameter $\gamma$ based on practical needs under the guidance of the above theoretical results. {\color{blue}For example, for safety-critical systems prioritizing rapid recovery, one may select an admissible $\gamma$ and the smallest fixed $\omega$ satisfying \eqref{Theorem_1_result_4} for the prescribed recovery horizon and violation bound. For performance-oriented systems prioritizing convergence accuracy, a decaying buffer is preferred.}
	
	\section{DRAPs with a Global Equality Constraint}
	\label{sec:extension}
	In this section, we will extend our proposed algorithm to DRAPs with a global equality constraint and also theoretically prove its convergence properties. {\color{blue}Although such a constraint can be equivalently represented by two inequalities, directly exploiting the equality structure yields a simpler specialization of DanyRA and enables explicit establishment of equality preservation and linear convergence to the exact optimal solution.}
	
	Consider the DRAPs with a global equality constraint
	\begin{align}
		\begin{split}
			& \min_{\mathbf{x}} \quad  f(\mathbf{x})=\sum\limits_{i=1}^n \ f_i(x_i) \\
			&  \ \ \text{s.t.} \quad \sum_{i=1}^n A_ix_i=\sum_{i=1}^nd_i.
		\end{split}\label{DRAP}
	\end{align}
	Similarly, by introducing an auxiliary variable $\mathbf{y}$, the DRAP (\ref{DRAP}) can be  equivalently written as 
	\begin{align}
		\begin{split}
			& \min_{x_i} \quad  \sum\limits_{i=1}^n \ f_i(x_i)  \\
			&  \ \text{s.t.} \ \ A_ix_i\!-\!d_i\!+\!\sum\nolimits_{j\in\mathcal{N}_i}w_{i,j}(y_i\!-\!y_j)\!=\!0, \ \forall i\in \mathcal{N},
		\end{split} \label{DRAP_1}
	\end{align}
	which is the special case of (\ref{IDRAP_1}) with known $\delta^\star=0$. Similarly, the equality version of DanyRA (Eq-DanyRA) algorithm is derived as summarized in Algorithm \ref{alg1}, where the set $
	\mathcal{K}^{\prime x}_i(x_{i,k})=\{x_i\in\mathbb{R}^p| A_i x_i=A_i x_{i,k} -\gamma \big(A_ix_i-d_i+\sum\nolimits_{j\in\mathcal{N}_i}w_{i,j}(y_{i,k+1}\!-\!y_{j,k+1})\big)\}$.
	\begin{algorithm}
		\renewcommand{\algorithmicrequire}{\textbf{Input:}}
		\footnotesize
		\caption{\bf : Eq-DanyRA algorithm}
		\label{alg1}
		\begin{algorithmic}[1]
			\STATE {Parameters: 
				proper stepsizes~$\alpha,\beta,\eta,\gamma>0$;\\
				Initialization: Arbitrary $x^\prime_{i,0},x_{i,0}\in \mathbb{R}^{p},\lambda_{i,0}\in \mathbb{R}^{m}$, and $y_{i,0}=0$.}
			\FOR{$k=0, 1, \cdots, K-1 $}
			\STATE \emph{Information exchange between agents:}  \\
			$\bar{\lambda}_{i,k}\!=\!\sum\limits_{j\in\mathcal{N}_i}w_{i,j}(\lambda_{i,k}\!-\!\lambda_{j,k})$, $\bar{y}_{i,k}=\sum\limits_{j\in\mathcal{N}_i}w_{i,j}(y_{i,k}-y_{j,k}),$\\ 
			$\bar{z}_{i,k}=\sum\limits_{j\in\mathcal{N}_i}w_{i,j}((z_{i,k}-d_i)-(z_{j,k}-d_j))$,\\
			\STATE \emph{Intermediate variable update:}  \\
			$x^\prime_{i,k+1}=x^\prime_{i,k}-\alpha(\nabla f_i(x^\prime_{{i,k}})+A_i^\intercal(z_{i,k}-d_i+\lambda_{i,k}))$,\\
			\STATE \emph{Auxiliary variable update:} \\
			$ y_{i,k+1}=y_{i,k}-\alpha(\bar{z}_{i,k}+\bar{\lambda}_{i,k})$,\\
			\STATE \emph{Dual variable update:} \\
			$\bar{y}_{i,k+1}=\sum\limits_{j\in\mathcal{N}_i}w_{i,j}(y_{i,k+1}-y_{j,k+1}),$\\
			$z_{i,k+1}=A_ix^\prime_{i,k+1}+\bar{y}_{i,k+1}$,\\
			$\lambda_{i,k+1}=\lambda_{i,k}+\beta\big(z_{i,k+1}-d_i-\eta A_i(A_i^\intercal\lambda_{i,k}+\nabla f_i(x^\prime_{{i,k}}))\big)$,\\	
			\STATE \emph{Decision variable update:}  \\
			$ x_{i,k+1}=\arg\min_{x_i\in\mathcal{K}^{\prime x}_i(x_{i,k})} \frac{1}{2}\Vert x_i-x^\prime_{i,k+1}\Vert^2.$
			\ENDFOR
		\end{algorithmic}
	\end{algorithm} 
	
	{\color{blue}For this equality-constrained specialization, both the preservation of the global equality constraint and the linear convergence rate can be established.} To do this, similarly, it is claimed that an optimal solution pair of (\ref{DRAP_1}) is attained if and only if the following conditions hold:
	\begin{align}
		& \nabla\! f(\mathbf{x}^\star)\!+\!\mathbf{A}^\intercal\!\boldsymbol{\lambda}^\star\!=\!0, \ \boldsymbol{\mathcal{L}}\boldsymbol{\lambda}^\star\!=\!0, \ \mathbf{A}\!\mathbf{x}^\star\!+\!\boldsymbol{\mathcal{L}}\mathbf{y}^\star\!-\!\mathbf{d}\!=\!0. \label{KKT_equality}
	\end{align}
	Then, it is given in Proposition 4 that any fixed point of Eq-DanyRA algorithm is consistent with the optimal solution of DRAP (\ref{DRAP_1}). This result can be easily obtained from  Proposition 2 and the detailed proof is omitted. 
	
	\emph{Proposition 4}: Under Assumption \ref{connected}, any fixed point of the
	proposed Eq-DanyRA algorithm, denoted as $(\hat{\mathbf{x}},\hat{\mathbf{y}},\hat{\boldsymbol{\lambda}})$ satisfies the optimality conditions
	given in (\ref{KKT_equality}).
	
	On this basis, we will establish convergence rate of Eq-DanyRA algorithm based on a modified metric
	$V_{k}=\Vert \tilde{\mathbf{x}}^\prime_{k} \Vert^2+\Vert\tilde{\mathbf{y}}_{k} \Vert^2+\frac{\alpha}{\beta}\Vert \tilde{\boldsymbol{\lambda}}_{k} \Vert^2 +\frac{\alpha(1-3\beta)}{2}\Vert\tilde{\mathbf{z}}_k\Vert^2+\frac{\underline{\sigma}_\mathbf{A}^2\alpha(1-3\beta)}{8\gamma^2}\Vert\mathbf{x}_{k}\!-\!\mathbf{x}^\prime_{k}\Vert^2$, where ${\mathbf{z}}_k=\mathbf{A}\mathbf{x}^\prime_k+\boldsymbol{\mathcal{L}}\mathbf{y}_k$ in this case. 
	
	{\emph{Theorem 3:}  Under Assumptions \ref{connected}-\ref{strong_convexity} and conditions in (\ref{conditions_Theorem_1}), if the following condition is satisfied:
		\begin{align}		
			\alpha< \frac{8\mu\!-\!4\eta\ell^2(3\beta\eta\!+\!1)+(1-3\beta)\bar{\sigma}_\mathbf{A}^2}{8\ell^2}, 	\label{conditions_Lemma_2}
		\end{align}
		the proposed Eq-DanyRA algorithm satisfies:
		
		\romannumeral1) $\{\mathbf{x}_{k}\}_{k\geq0}$ converges linearly to the optimal solution of DRAP (\ref{DRAP}) as 
		\begin{align*}
			\Vert\mathbf{x}_{k+1}\!-\!\hat{\mathbf{x}}\Vert^2 \leq\frac{16\gamma^2V_0{\theta^\prime}^k}{\min\{8\gamma^2,{\underline{\sigma}_\mathbf{A}^2\alpha(1-3\beta)}\}},
		\end{align*}
		where the constant  $\theta^\prime\!=\!\max\{1\!+\!\alpha(2\ell^2\alpha\!-\!2\mu\!+\!\eta\ell^2(3\beta\eta\!+\!1)+\frac{(3\beta-1)\bar{\sigma}_\mathbf{A}^2}{4}),\frac{8\!+\!\alpha(3\beta-1)\underline{\sigma}_\mathcal{L}^2}{8},1\!+\!\beta\eta(3\beta\eta\bar{\sigma}_\mathbf{A}^2\!-\!\underline{\sigma}_\mathbf{A}^2),\frac{1}{2},{4(1\!-\!\gamma)^2\kappa_\mathbf{A}^2}\}$ is strictly less than one.
		
		\romannumeral2) when the coupled constraint is violated for some $k_0\geq 1$, the violation $C^{vio}_{k_0}=\Vert\sum_{i=1}^n (A_ix_{i,k_0}-d_{i})\Vert_\infty$ decreases linearly to zero, i.e.,
		\begin{align}
			\sum_{i=1}^n \left(A_ix_{i,k_0+t}-d_{i}\right) \preccurlyeq (1-\gamma)^t C^{vio}_{k_0}\mathbf{1}, \ \forall t\geq 1;
		\end{align}
		\romannumeral3) once the coupled constraint is satisfied, it remains satisfied thereafter, i.e., if $\sum_{i=1}^n A_ix_{i,k_0}=\sum_{i=1}^nd_{i}$ for some $k_0\geq 1$, then 
		\begin{align}
			\sum_{i=1}^n A_ix_{i,k}=\sum_{i=1}^nd_{i}, \ \forall k\geq k_0.
		\end{align}
		\emph{Proof:} See Appendix \ref{Proof_of_Theorem_3}.
		$\hfill \blacksquare$}
	
	{\color{blue}By exploiting the equality structure, we strengthen the Lyapunov descent estimate for the inequality-constrained case to a contraction for Eq-DanyRA, yielding an explicit linear convergence rate. The coupled equality constraint is preserved once satisfied, and any violation decays geometrically between disruptions.}
	However, robustness to constraint violation cannot be ensured for Eq-DanyRA algorithm, which requires further investigation in future work.

	\section{Numerical Experiments}
	\label{sec:numerical_experiments}
	In this section, we conduct numerical experiments to verify our theoretical results. We consider the DRAP for
	industrial Internet of Things based control systems \cite{2024_TII_Distributed_Multidomain_Resource_Allocation} consisting of $n=14$ real-time computation tasks awaiting completion. To ensure all of them being completed before deadlines, computation resources require proper allocation under the upper bound $r^{\max}$, and the schedulability condition (see Theorem 4.2 in \cite{Hard_Real_Time_Computing_Systems}) must be satisfied. Moreover, taken the total computational energy consumption of all tasks as the cost function, we have that this problem can be rewritten as DRAP (\ref{IDRAP}), where $f_i=x_i^\intercal P_ix_i-Q_i^\intercal x_i$, $A_i=blkdiag(1,C_i)$, $\sum_{i=1}^n d_i=[r^{\max} \ 1]^\intercal$, $x_i=[r_i \ \frac{1}{t_i}]^\intercal$, $r_i$ is the computation resources allocated to task $i$ and $t_i$ denotes its period. Moreover, the coefficients $C_i>0$, and $P_i,Q_i\succ0$ are randomly generated.
	Regarding the communication network, we generate an undirected connected graph by adding random links to a ring network.  Experiments are implemented on a computer with Intel-i7 2.1 GHz CPU and 16 GB of RAM.
	
\begin{figure}[thpb]
		\centering
		\includegraphics[width=6.5cm]{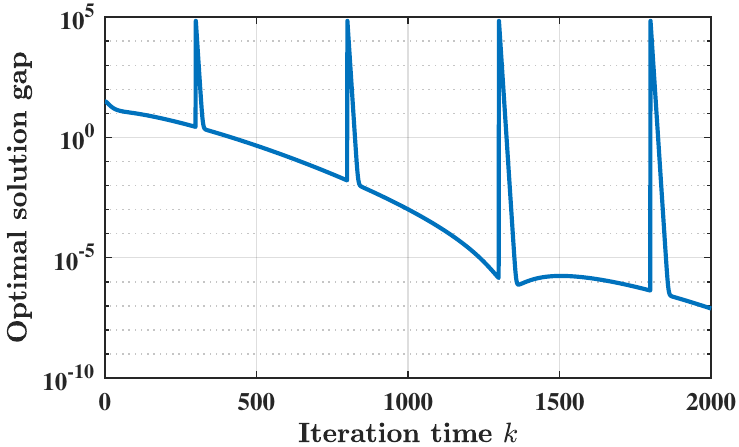}
		\caption{\color{blue}Evolutions of optimality gap with respect to iteration time $k$.}
		\label{Gap_In_DanyRA}
	\end{figure}

\begin{figure}[thpb]
		\centering
		\includegraphics[width=6.5cm]{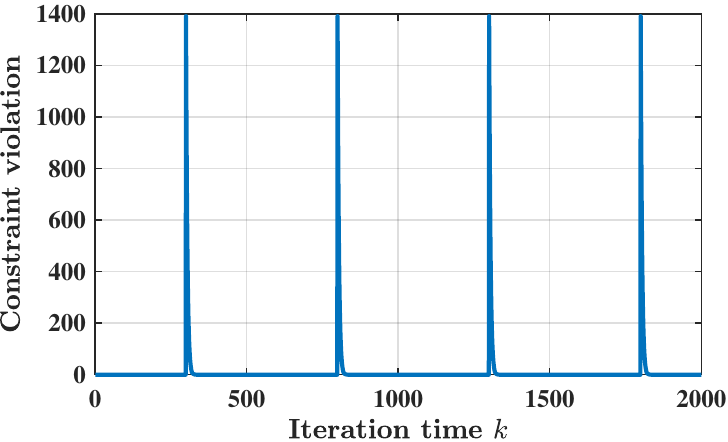}
		\caption{\color{blue}Evolutions of constraint violation with respect to iteration time $k$.}
		\label{Violation_In_DanyRA}
	\end{figure}

	Set the upper bound of the computation resources $r^{\max}=70$, and {\color{blue}the local budget contribution $d_i=[5 \ \frac{1}{14}]^\intercal,i=1,\cdots,14$}. We first verify the effectiveness of DanyRA algorithm. Choose the stepsizes $\alpha=0.01$, $\beta=0.02$, $\gamma=0.2$, $\eta=0.1$, and $\omega=0$, with $x_{i,0}=d_i$, for all agents. {\color{blue}Every $500$ iterations, we perturb each decision variable by $[50 \ 50]^\intercal$, causing repeated violations of the global constraint.}  Fig. \ref{Gap_In_DanyRA} plots the historical evolution of the optimality gap $\Vert\mathbf{x}_k-\mathbf{x}^\star\Vert^2$ of DanyRA algorithm. The figures show that the decision $\mathbf{x}$ converges fast to the optimal solution of DRAPs (\ref{IDRAP}) in the absence of interruption. Fig. \ref{Violation_In_DanyRA} plots the historical evolution of the constraint violation quantified by the 1-norm $||\max\{\sum_{i=1}^n(A_ix_i-d_i),0\}||_1$. It is shown that once the derived solution is feasible, its feasibility can be maintained thereafter; even if a constraint violation occurs, it decays rapidly to zero.  {\color{blue}These results illustrate repeated feasibility recovery when the inter-disruption intervals are sufficiently long.}
	
\begin{figure}[thpb]
		\centering
		\includegraphics[width=6.5cm]{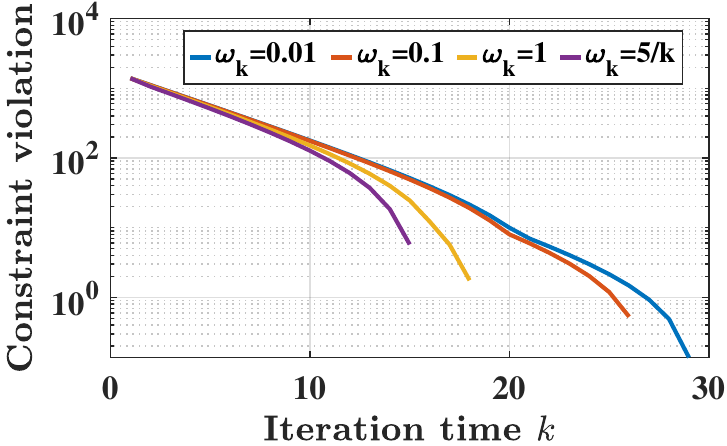}
		\caption{Evolutions of constraint violation with respect to iteration time $k$.}
		\label{Violation_In_DanyRA_diff_Delta}
	\end{figure}

\begin{figure}[thpb]
		\centering
		\includegraphics[width=6.5cm]{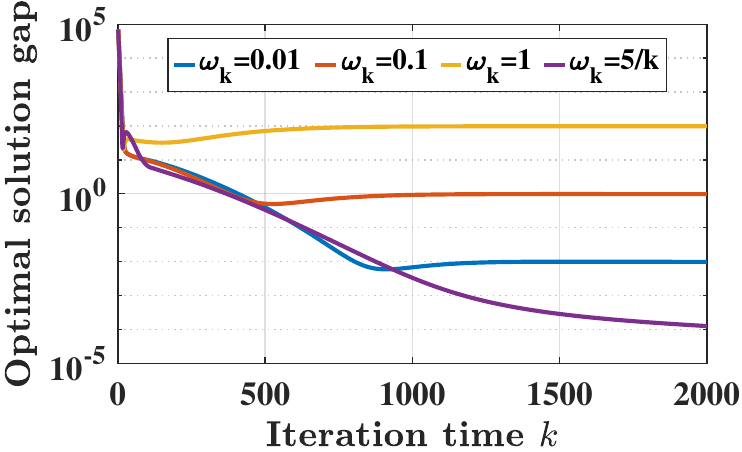}
		\caption{Evolutions of optimality gap with respect to iteration time $k$.}
		\label{Gap_In_DanyRA_diff_Delta}
	\end{figure}

	In the case with nonzero $\omega_k$, we consider the following choices: $\omega_k=10^{-2}$, $10^{-1}$, $1$, $5/k$. Here, we initialize  $x_{i,0}=d_i+[50 \ 50]^\intercal,\forall i\in\mathcal{N}$ and the corresponding results of the experiments are given as follows. It can be seen from Fig. \ref{Violation_In_DanyRA_diff_Delta} that the constraint violation can be eliminated within finite time. Moreover, combining this result with the one shown in Fig. \ref{Gap_In_DanyRA_diff_Delta} that when $\omega_k$ is set to a positive constant, a larger 
	$\omega_k$ reduces the required number of iterations for feasibility recovery but increases the corresponding optimality gap. Thus, the results in Theorem 1 and Theorem 2 are numerically verified. Additionally, by selecting an appropriate decaying buffer, it can be observed that both convergence accuracy and rapid feasibility recovery can be achieved simultaneously.
	
\begin{figure}[thpb]
		\centering
		\includegraphics[width=6.5cm]{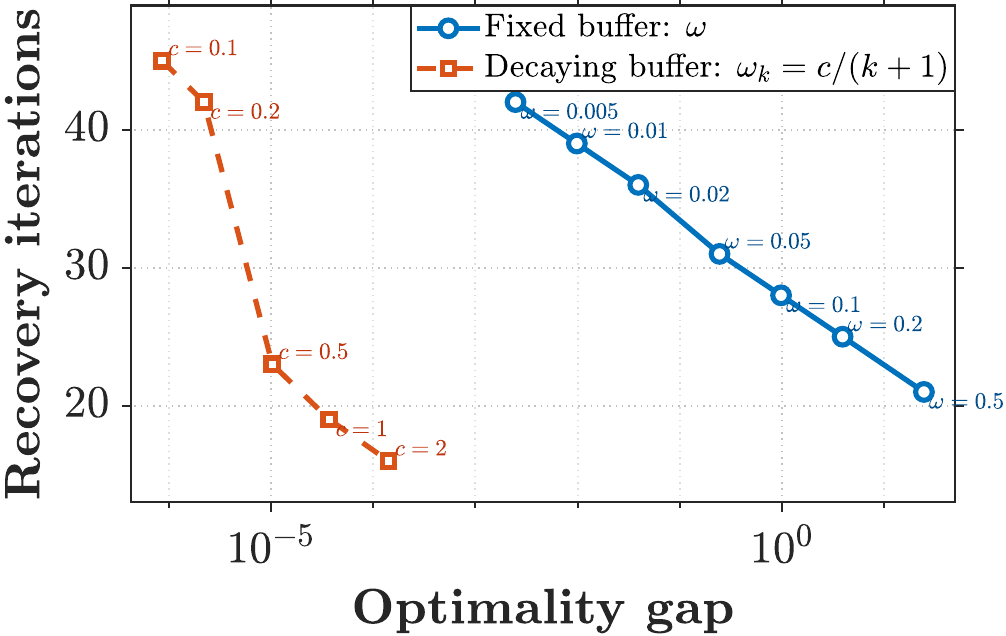}
		\caption{\color{blue}Accuracy--robustness trade-off for different buffers.}
		\label{fig:accuracy_robustness_tradeoff}
	\end{figure}

	{\color{blue}To quantify this trade-off, we repeat the experiment by varying only the fixed buffer $\omega$ or the coefficient $c$ in the decaying buffer $\omega_k=c/(k+1)$.  Convergence accuracy is measured by the optimality gap, while violation robustness is measured by the recovery time required for the violation to remain below $10^{-10}$. Fig.~\ref{fig:accuracy_robustness_tradeoff} shows that a larger fixed buffer accelerates recovery at the cost of accuracy, whereas decaying buffers achieve comparable recovery with substantially smaller optimality gaps. Moreover, the observed reduction in recovery time with a larger fixed buffer supports its use in safety-critical systems prioritizing rapid feasibility recovery. The smaller optimality gaps achieved with decaying buffers make them suitable for performance-oriented systems prioritizing convergence accuracy.}
	
\begin{figure}[thpb]
		\centering
		\includegraphics[width=6.5cm]{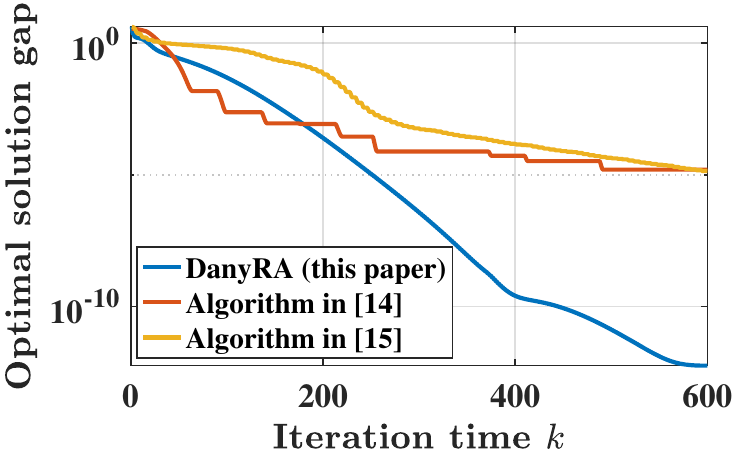}
		\caption{Evolutions of optimality gap with respect to iteration time $k$.}
		\label{Algorithms_comparison}
	\end{figure}

	Next, we compare the proposed DanyRA algorithm with the state-of-the-art algorithms in \cite{2024_arxiv_Achieving_violation_free_distributed_optimization_under_coupling_constraints,2023_CDC_Distributed_and_Anytime_Algorithm_for_Network_Optimization_Problems_with_Separable_Structure}, in terms of the overall optimality gap $\Vert\mathbf{x}_k-\mathbf{x}^\star\Vert^2$. {\color{blue}Each algorithm is independently tuned through a coarse-to-fine search over its admissible parameter ranges, and the best-observed stable configuration is used. To be specific, we set $\alpha=3.67*10^{-3}$, $\beta=0.67$, $\gamma=0.69$, $\eta=0.074$ and $\omega_k=0, \forall k\geq0$ for DanyRA algorithm; $\gamma=1.32*10^{-6}$ for Algorithm in \cite{2024_arxiv_Achieving_violation_free_distributed_optimization_under_coupling_constraints}; $\tau=0.27$, $\beta=5.54$, and $\epsilon=0.15$ for Algorithm in \cite{2023_CDC_Distributed_and_Anytime_Algorithm_for_Network_Optimization_Problems_with_Separable_Structure}.}  As shown in Fig. \ref{Algorithms_comparison}, the optimality gap of the proposed DanyRA algorithm decreases faster than that of the other two algorithms, verifying its effectiveness in achieving optimality for DRAPs.
	
	{\color{blue}Finally, we compare the computation time of the algorithms. {For this comparison, all three algorithms are independently tuned at each network size.} {The stopping criterion is satisfied when both the optimality gap and the constraint violation remain below} \({10^{-5}}\) and \(10^{-10}\), respectively, for 20 consecutive iterations. Each experiment is repeated fifteen times, and the median computation time is reported in Table \ref{computation_time}. {As shown in Table \ref{computation_time}, the proposed DanyRA algorithm achieves the lowest median computation time among the three algorithms at every tested network size under the same stopping criterion.}} {This verifies that the proposed algorithm has lower computational complexity.}
	
	\begin{table}[!ht]
		\renewcommand{\arraystretch}{1.3}
		\centering
		\caption{\color{blue}Computation time of algorithms}
		\label{computation_time}
		\resizebox{\linewidth}{!}{\begin{tabular}{|c|c|c|c|c|c|c|}
				\hline
				\multicolumn{2}{|c|}{Number of agents} & 20 & 30 & 40 & 50   \\
				\hline
				\multirow{3}{*}{Computation} & DanyRA algorithm  & {3.684} & {9.600} & {13.751} & {19.077}    \\
				\cline{2-6}  
				& Algorithm in \cite{2024_arxiv_Achieving_violation_free_distributed_optimization_under_coupling_constraints}  & {16.455} & {24.840} & {35.129} & {46.886}  \\
				\cline{2-6}
				time (s) & Algorithm in \cite{2023_CDC_Distributed_and_Anytime_Algorithm_for_Network_Optimization_Problems_with_Separable_Structure}  & {10.921} & {11.652} & {15.706} & {36.864} \\
				\hline
		\end{tabular}}
	\end{table}
	
	\section{Conclusion}
	\label{sec:sonclusion}
	In this paper, by synthesizing a safe feedback controller using control barrier functions, a distributed anytime-feasible resource allocation (DanyRA) algorithm has been proposed. It achieves optimality for DRAPs with a coupled inequality constraint while ensuring anytime feasibility. Under standard assumptions, its exact convergence to the optimal solution has been proved. Considering potential constraint violation, the algorithm has been redesigned to achieve robustness for maintaining or recovering feasibility before a pre-defined deadline. Then, the trade-off between convergence accuracy and violation robustness has been characterized. The proposed DanyRA algorithm has been further extended to DRAPs with a coupled equality constraint, and its anytime feasibility and linear convergence rate have been proved. Finally, numerical experiments have validated the theoretical results. {\color{blue}Future work will focus on extending the recovery analysis to persistent or high-frequency disturbances, developing adaptive buffer mechanisms, and generalizing DanyRA to directed or time-varying networks and more general objective functions.}

	\newcounter{saveeqn}
	\setcounter{saveeqn}{\value{equation}}
	
	\appendix
	\renewcommand{\thesection}{\Alph{section}}
	\setcounter{section}{0}
	
	\setcounter{equation}{\value{saveeqn}}
	\counterwithout{equation}{section} 
	
	\section{Proof of Proposition 1}
	\label{Proof_of_Proposition_1}
	On the one hand, since $\sum_{i\in\mathcal{N}}\sum_{j\in\mathcal{N}_i}w_{i,j}(y_i-y_j)=0$ under Assumption \ref{connected} and $\boldsymbol{\delta}\succcurlyeq0$, it can be easily verified that for any feasible solution $(\mathbf{x},\mathbf{y})$ of problem (\ref{IDRAP_1}), $\mathbf{x}$ is a feasible solution to (\ref{IDRAP}). On the other hand, for any feasible solution $\mathbf{x}$ of problem (\ref{IDRAP}), we have $\hat{\mathbf{1}}^\intercal(\mathbf{A}\mathbf{x}-\mathbf{d})\preccurlyeq0$. Thus, $\hat{\mathbf{1}}^\intercal(\mathbf{A}\mathbf{x}+\boldsymbol{\delta}-\mathbf{d})=0$ holds for some $\boldsymbol{\delta}\succcurlyeq0$, which means $\mathbf{A}\mathbf{x}+\boldsymbol{\delta}-\mathbf{d}$ is in the range space of $\boldsymbol{\mathcal{L}}$. Thus, there exists some $\mathbf{y}$ such that $(\mathbf{x},\mathbf{y})$ is feasible to (\ref{IDRAP_1}). Combining these results with the fact that the problems (\ref{IDRAP}) and (\ref{IDRAP_1}) have the same cost function $f(\mathbf{x})=\sum\limits_{i=1}^n \ f_i(x_i)$ completes the proof. $\hfill\blacksquare$
	
	\section{Proof of Proposition 2}
	\label{Proof_of_Proposition_2}
	It is not difficult to check that any fixed point
	of the DanyRA algorithm satisfies
	\begin{align}
		&  \mathbf{A}^\intercal(\mathbf{A}\hat{\mathbf{x}}^\prime+\boldsymbol{\mathcal{L}}\hat{\mathbf{y}}+\hat{\boldsymbol{\delta}}-\mathbf{d})=-\left(\nabla f(\hat{\mathbf{x}}^\prime)+\mathbf{A}^\intercal\hat{\boldsymbol{\lambda}}\right), \label{fix_In_DanyRA_1}\\
		& \boldsymbol{\mathcal{L}}\left(\mathbf{A}\hat{\mathbf{x}}^\prime+\boldsymbol{\mathcal{L}}\hat{\mathbf{y}}+\hat{\boldsymbol{\delta}}-\mathbf{d}\right)+\boldsymbol{\mathcal{L}}\hat{\boldsymbol{\lambda}}=0,\label{fix_In_DanyRA_2}\\
		&\max\big\{\hat{\boldsymbol{\delta}}-\alpha(\mathbf{A}\hat{\mathbf{x}}^\prime+\boldsymbol{\mathcal{L}}\hat{\mathbf{y}}+\hat{\boldsymbol{\delta}}-\mathbf{d}+\hat{\boldsymbol{\lambda}}),0 \big\}=\hat{\boldsymbol{\delta}}, \label{fix_In_DanyRA_3}\\
		&\mathbf{A}\hat{\mathbf{x}}^\prime+\boldsymbol{\mathcal{L}}\hat{\mathbf{y}}+\hat{\boldsymbol{\delta}}-\mathbf{d}=\eta\mathbf{A}\left(\nabla f(\hat{\mathbf{x}}^\prime)+\mathbf{A}^\intercal\hat{\boldsymbol{\lambda}}\right), \label{fix_In_DanyRA_4} \\
		&\mathbf{A}\hat{\mathbf{x}}+\boldsymbol{\mathcal{L}}\hat{\mathbf{y}}+\hat{\boldsymbol{\delta}}-\mathbf{d}=0. \label{fix_In_DanyRA_5}
	\end{align} 
	Combining (\ref{fix_In_DanyRA_1}) with (\ref{fix_In_DanyRA_4}) yields $\mathbf{A}\hat{\mathbf{x}}^\prime+\boldsymbol{\mathcal{L}}\hat{\mathbf{y}}+\hat{\boldsymbol{\delta}}-\mathbf{d}=\nabla f(\hat{\mathbf{x}}^\prime)+\mathbf{A}^\intercal\hat{\boldsymbol{\lambda}}=0$. Then, substituting this relation into relations (\ref{fix_In_DanyRA_2}) and (\ref{fix_In_DanyRA_5}) yields $\boldsymbol{\mathcal{L}}\hat{\boldsymbol{\lambda}}=0$ and $\mathbf{A}\hat{\mathbf{x}}^\prime=\mathbf{A}\hat{\mathbf{x}}$, respectively. From $\mathbf{A}\hat{\mathbf{x}}^\prime=\mathbf{A}\hat{\mathbf{x}}$, we have that $\hat{\mathbf{x}}^\prime$ is also a feasible solution of the optimization problem in step $7$ and thus $\hat{\mathbf{x}}^\prime=\hat{\mathbf{x}}$ holds due to its objective function. 
	
	It can be derived from (\ref{fix_In_DanyRA_3}) that $\max\{\hat{\boldsymbol{\delta}}-\alpha\hat{\boldsymbol{\lambda}},0\}=\hat{\boldsymbol{\delta}}$. This implies that if $\hat{\boldsymbol{\delta}}=0$, then $\hat{\boldsymbol{\lambda}}\geq 0$; if $\hat{\boldsymbol{\delta}}>0$, then $\hat{\boldsymbol{\lambda}}=0$. Thus, we have $\hat{\boldsymbol{\lambda}}\geq0$, $\hat{\boldsymbol{\delta}}\geq 0$, and ${\hat{\boldsymbol{\lambda}}}^\intercal\hat{\boldsymbol{\delta}}=0$. 
	
	To sum, $(\hat{\mathbf{x}},\hat{\mathbf{y}},\hat{\boldsymbol{\delta}},\hat{\boldsymbol{\lambda}})$ satisfies the optimality condition (\ref{KKT_inequality}). The proof is completed. $\hfill\blacksquare$

	\section{Proof of Proposition 3}
	\label{Proof_of_Proposition_3}
	According to Proposition 2, we can set a fixed point of the DanyRA algorithm with $\omega=0$, to be an optimal solution of (\ref{IDRAP}). Thus, it can be derived that relations (\ref{KKT_inequality}) and $\boldsymbol{\lambda}^\star\succcurlyeq 0$ hold. 
	
	It can be derived from (\ref{fix_In_DanyRA_3}) with $\omega>0$ that $\max\big\{\hat{\boldsymbol{\delta}}-\hat{\omega}\mathbf{1}-\alpha\hat{\boldsymbol{\lambda}},0\big\}=\hat{\boldsymbol{\delta}}-\hat{\omega}\mathbf{1}$, where $\hat{\omega}$ is the limit point of $\omega$ ($\omega_{k}$). Similar to the proof of Proposition 2, we have
	\begin{align}
		& \nabla f(\hat{\mathbf{x}})+\mathbf{A}^\intercal\hat{\boldsymbol{\lambda}}=0, \ \boldsymbol{\mathcal{L}}\hat{\boldsymbol{\lambda}}=0,  \label{Proposition_3_1} \\ 
		& \hat{\boldsymbol{\lambda}}\succcurlyeq 0, \  \hat{\boldsymbol{\delta}}\geq\hat{\omega}\mathbf{1}, \ {\hat{\boldsymbol{\lambda}}}^\intercal(\hat{\boldsymbol{\delta}}-\hat{\omega}\mathbf{1})=0, \label{Proposition_3_2} \\ &\mathbf{A}\hat{\mathbf{x}}+\boldsymbol{\mathcal{L}}\hat{\mathbf{y}}+\hat{\boldsymbol{\delta}}-\mathbf{d}=0. \label{Proposition_3_3} 
	\end{align}
	From (\ref{KKT_inequality_1}),  (\ref{Proposition_3_1}) and Assumption \ref{strong_convexity}, we have
	\begin{align}
		\mu\Vert\hat{\mathbf{x}}\!-\!\mathbf{x}^\star\Vert^2\!\leq&(\nabla f(\hat{\mathbf{x}})-\nabla f(\mathbf{x}^\star))^\intercal(\hat{\mathbf{x}}-\mathbf{x}^\star) \notag \\
		=&-(\hat{\boldsymbol{\lambda}}-\boldsymbol{\lambda}^\star)^\intercal \mathbf{A}(\hat{\mathbf{x}}-\mathbf{x}^\star) \notag \\
		=& (\hat{\boldsymbol{\lambda}}\!-\!\boldsymbol{\lambda}^\star)^\intercal\hat{\omega}\mathbf{1} \!+\! (\hat{\boldsymbol{\lambda}}\!-\!\boldsymbol{\lambda}^\star)^\intercal (\hat{\boldsymbol{\delta}}\!-\!\hat{\omega}\mathbf{1}\!-\!\boldsymbol{\delta}^\star), \label{Proposition_3_4} 
	\end{align}
	where the last equality uses relations (\ref{KKT_inequality_2}) and (\ref{Proposition_3_3}).
	
	Considering the last term in (\ref{Proposition_3_4}), using relations (\ref{KKT_inequality_2}) and (\ref{Proposition_3_2}), it can be derived that
	\begin{align}
		&(\hat{\boldsymbol{\lambda}}\!-\!\boldsymbol{\lambda}^\star)^\intercal (\hat{\boldsymbol{\delta}}\!-\!\hat{\omega}\mathbf{1}\!-\!\boldsymbol{\delta}^\star) \notag \\
		=&\hat{\boldsymbol{\lambda}}^\intercal(\hat{\boldsymbol{\delta}}\!-\!\hat{\omega}\mathbf{1})+{\boldsymbol{\lambda}^\star}^\intercal\boldsymbol{\delta}^\star\!-\!{\boldsymbol{\lambda}^\star}^\intercal(\hat{\boldsymbol{\delta}}\!-\!\hat{\omega}\mathbf{1})-\hat{\boldsymbol{\lambda}}^\intercal\boldsymbol{\delta}^\star\leq0. \label{Proposition_3_5}
	\end{align}
	Substituting (\ref{Proposition_3_5}) into (\ref{Proposition_3_4}) and using (\ref{KKT_inequality_1}),  (\ref{Proposition_3_1}) yield $
	\mu\Vert\hat{\mathbf{x}}\!-\!\mathbf{x}^\star\Vert^2\!\leq \Vert\hat{\boldsymbol{\lambda}}-\boldsymbol{\lambda}^\star\Vert\Vert\hat{\omega}\mathbf{1}\Vert \leq \frac{\sqrt{n}\hat{\omega}}{\underline{\sigma}_\mathbf{A}}\Vert \mathbf{A}^\intercal\hat{\boldsymbol{\lambda}}-\mathbf{A}^\intercal\boldsymbol{\lambda}^\star\Vert=\frac{\sqrt{n}\hat{\omega}}{\underline{\sigma}_\mathbf{A}}\Vert\nabla f(\hat{\mathbf{x}})-\nabla f(\mathbf{x}^\star)\Vert\leq\frac{\ell\sqrt{n}\hat{\omega}}{\underline{\sigma}_\mathbf{A}}\Vert\hat{\mathbf{x}}\!-\!\mathbf{x}^\star\Vert$. Combining this result with the fact that $\hat{\omega}=\omega$ completes the proof. $\hfill\blacksquare$
	
	\section{Proof of Lemma 1}
	\label{Proof_of_Lemma_1}
	For conciseness, we first rewrite the DanyRA algorithm into a compact form as 
	\begin{align}
		& \mathbf{x}^\prime_{k+1}=\mathbf{x}^\prime_{k}-\alpha(\nabla f(\mathbf{x}^\prime_k)+\mathbf{A}^\intercal(\mathbf{z}_k-\mathbf{d})+\mathbf{A}^\intercal\boldsymbol{\lambda}_k), \label{compact_In_DanyRA_1} \\
		& \mathbf{y}_{k+1}=\mathbf{y}_{k}-\alpha \boldsymbol{\mathcal{L}}(\mathbf{z}_k-\mathbf{d}+\boldsymbol{\lambda}_k), \label{compact_In_DanyRA_2} \\
		& \boldsymbol{\delta}_{k+1}=\max\{\boldsymbol{\delta}_k-\alpha(\mathbf{z}_k-\mathbf{d}+\boldsymbol{\lambda}_k),{\omega}\}, \label{compact_In_DanyRA_3} \\
		& \boldsymbol{\lambda}_{k+1}=\!\boldsymbol{\lambda}_k\!+\!\beta\left(\mathbf{z}_{k+1}\!-\!\mathbf{d}\!-\!\eta\mathbf{A}(\mathbf{A}^\intercal\boldsymbol{\lambda}_k\!+\!\nabla f(\mathbf{x}^\prime_k))\right), \label{compact_In_DanyRA_4} \\
		& \mathbf{x}_{k+1}=\arg\min\nolimits_{\mathbf{x}\in\mathcal{K}^x(\mathbf{x}_{k})} \frac{1}{2}\Vert \mathbf{x}-\mathbf{x}^\prime_{k+1}\Vert^2,\label{compact_In_DanyRA_5}
	\end{align}
	where $\mathcal{K}^x(\mathbf{x}_{k})=\{\mathbf{x}\in\mathbb{R}^{np}|\mathbf{A}\mathbf{x}=\mathbf{A}\mathbf{x}_{k}-\gamma(\mathbf{A}\mathbf{x}_k+\boldsymbol{\mathcal{L}}\mathbf{y}_{k+1}+\boldsymbol{\delta}_{k+1}-\mathbf{d})+(1-\gamma)(\boldsymbol{\delta}_k-\boldsymbol{\delta}_{k+1})\}$.
	
	It can be derived from (\ref{compact_In_DanyRA_1}), (\ref{compact_In_DanyRA_2}) and (\ref{compact_In_DanyRA_4}) that
	\begin{align}
		\Vert\tilde{\mathbf{x}}^\prime_{k+1}\Vert^2 
		=&\Vert\tilde{\mathbf{x}}^\prime_{k}\Vert^2-\Vert\Delta{\mathbf{x}}^\prime_{k+1}\Vert^2-2\alpha\tilde{\mathbf{x}}_{k+1}^{\prime\intercal}\mathbf{A}^\intercal\tilde{\boldsymbol{\lambda}}_k \notag \\
		-&2\alpha\tilde{\mathbf{x}}_{k+1}^{\prime\intercal}(\nabla f(\mathbf{x}^\prime_k)-\nabla f(\hat{\mathbf{x}}^\prime)+\mathbf{A}^\intercal\tilde{\mathbf{z}}_k),  \label{proof_In_DanyRA_1} \\
		\Vert\tilde{\mathbf{y}}_{k+1}\Vert^2 
		=&\Vert\tilde{\mathbf{y}}_{k}\Vert^2-\Vert\Delta{\mathbf{y}}_{k+1}\Vert^2-2\alpha\tilde{\mathbf{y}}_{k+1}^\intercal\tilde{\boldsymbol{\lambda}}_k \notag \\
		-&2\alpha\tilde{\mathbf{y}}_{k+1}^\intercal\boldsymbol{\mathcal{L}}(\tilde{\mathbf{z}}_k+\tilde{\boldsymbol{\lambda}}_k), \label{proof_In_DanyRA_2} \\
		\Vert\tilde{\boldsymbol{\lambda}}_{k+1}\Vert^2 =&\Vert\tilde{\boldsymbol{\lambda}}_{k}\Vert^2+\Vert\Delta{\boldsymbol{\lambda}}_{k+1}\Vert^2+2\beta\tilde{\boldsymbol{\lambda}}_{k}^\intercal\tilde{\mathbf{z}}_{k+1} \notag\\
		-&2\beta\eta\tilde{\boldsymbol{\lambda}}_{k}^\intercal\mathbf{A}(\mathbf{A}^\intercal\tilde{\boldsymbol{\lambda}}_k+\nabla f(\mathbf{x}^\prime_k)-\nabla f(\hat{\mathbf{x}}^\prime)), \label{proof_In_DanyRA_3}
	\end{align}
	where we use relations (\ref{Proposition_3_1}) and (\ref{Proposition_3_3}).
	
	{\color{blue}As for the update of $\boldsymbol{\delta}$, it can be derived from (\ref{compact_In_DanyRA_3}) and the relation $\hat{\boldsymbol{\delta}}\succcurlyeq\omega\mathbf{1}$ that
		\begin{align*}
			\tilde{\boldsymbol{\delta}}_{k+1}^{\intercal}
			\left[
			\boldsymbol{\delta}_{k+1}
			-\big(\boldsymbol{\delta}_{k}
			-\alpha(\mathbf{z}_k-\mathbf{d}+\boldsymbol{\lambda}_k)\big)
			\right]
			\leq 0, 
		\end{align*}
		where the relation  $\left(\max\{\mathbf{x},{\omega}\}-\mathbf{y}\right)^{\intercal}\left(\max\{\mathbf{x},{\omega}\}-\mathbf{x}\right)\leq0, \forall \mathbf{x}\in\mathbb{R}^{nm}, \forall \mathbf{y}\in\{\mathbf{y}\vert \mathbf{y}-\omega\mathbf{1}\succcurlyeq0, \mathbf{y}\in\mathbb{R}^{nm} \}$ is used.
		
		From (\ref{compact_In_DanyRA_3}) and the definition of $\Delta\delta_{k+1}$, we have
		\begin{align*}
			\|\tilde{\boldsymbol{\delta}}_{k+1}\|^2-\|\tilde{\boldsymbol{\delta}}_k\|^2
			=
			-\|\Delta{\boldsymbol{\delta}}_{k+1}\|^2
			+
			2\tilde{\boldsymbol{\delta}}_{k+1}^{\top}\Delta{\boldsymbol{\delta}}_{k+1}.
		\end{align*}
		Combining the above two relations yields
		\begin{align}
			&\Vert\tilde{\boldsymbol{\delta}}_{k+1}\Vert^2-\Vert\tilde{\boldsymbol{\delta}}_{k}\Vert^2 \notag\\ 
			=&-\Vert\Delta{\boldsymbol{\delta}}_{k+1}\Vert^2-2\alpha \tilde{\boldsymbol{\delta}}_{k+1}^\intercal (\mathbf{z}_k-\mathbf{d}+\boldsymbol{\lambda}_k)\notag \\
			&+2\tilde{\boldsymbol{\delta}}_{k+1}^\intercal\left(\boldsymbol{\delta}_{k+1}-(\boldsymbol{\delta}_{k}-\alpha(\mathbf{z}_k-\mathbf{d}+\boldsymbol{\lambda}_k))\right) \notag \\
			\leq& -\Vert\Delta{\boldsymbol{\delta}}_{k+1}\Vert^2-2\alpha \tilde{\boldsymbol{\delta}}_{k+1}^\intercal (\tilde{\mathbf{z}}_k+\tilde{\boldsymbol{\lambda}}_k)-2\alpha \tilde{\boldsymbol{\delta}}_{k+1}^\intercal \hat{\boldsymbol{\lambda}}. \label{proof_In_DanyRA_4}
	\end{align}}
	The last term in the RHS of (\ref{proof_In_DanyRA_4}) can be rewritten as 
	\begin{align}
		-2\alpha \tilde{\boldsymbol{\delta}}_{k+1}^\intercal \hat{\boldsymbol{\lambda}}=-2\alpha{\boldsymbol{\delta}}_{k+1}^\intercal\hat{\boldsymbol{\lambda}}+2\alpha\hat{\boldsymbol{\delta}}^\intercal\hat{\boldsymbol{\lambda}}
		\leq 0, \label{proof_In_DanyRA_5}
	\end{align}
	where we use ${\boldsymbol{\delta}}_{k+1} \succcurlyeq {\omega}\mathbf{1}$, $\hat{\boldsymbol{\lambda}}\succcurlyeq 0$, and ${\hat{\boldsymbol{\lambda}}}^\intercal(\hat{\boldsymbol{\delta}}-{\omega}\mathbf{1})=0$. 
	
	Define $V_k^{\prime}=\Vert \tilde{\mathbf{x}}^\prime_{k} \Vert^2+\Vert\tilde{\mathbf{y}}_{k} \Vert^2+\Vert\tilde{\boldsymbol{\delta}}_{k} \Vert^2+\frac{\alpha}{\beta}\Vert \tilde{\boldsymbol{\lambda}}_{k} \Vert^2$. Combining relations (\ref{proof_In_DanyRA_1})-(\ref{proof_In_DanyRA_5}) yields
	\begin{align}
		&V_{k+1}^{\prime}-V_k^{\prime} 
		\leq  -\Vert\Delta{\mathbf{x}}^\prime_{k+1}\Vert^2-\Vert\Delta{\mathbf{y}}_{k+1}\Vert^2-\Vert\Delta{\boldsymbol{\delta}}_{k+1}\Vert^2 \notag\\
		&  -2\alpha\tilde{\mathbf{x}}_{k+1}^{\prime\intercal}(\nabla f(\mathbf{x}^\prime_k)\!-\!\nabla f(\hat{\mathbf{x}}^\prime))\!-\!2\alpha{\tilde{\mathbf{z}}_{k+1}}^\intercal\tilde{\mathbf{z}}_k  \!+\!\frac{\alpha}{\beta}\Vert{\boldsymbol{\lambda}}_{k+1}\!-\!{\boldsymbol{\lambda}}_{k}\Vert^2 \notag\\
		&-\!2\alpha\eta \Vert \mathbf{A}^\intercal \tilde{\boldsymbol{\lambda}}_k \Vert^2\!-2\alpha\eta\tilde{\boldsymbol{\lambda}}_{k}^\intercal\mathbf{A}(\nabla f(\mathbf{x}^\prime_k)-\nabla f(\hat{\mathbf{x}}^\prime)). \label{proof_In_DanyRA_6}
	\end{align}
	Then, we will give the upper bound of the RHS of (\ref{proof_In_DanyRA_6}) to prove the result (\ref{result_Lemma_1}). First, under Assumptions \ref{Lipschitz_smooth} and \ref{strong_convexity}, the second line in (\ref{proof_In_DanyRA_6}) can be upper bounded as 
	\begin{align}
		-&\!2\alpha\tilde{\mathbf{x}}_{k\!+\!1}^{\prime\intercal}(\nabla\! f(\mathbf{x}^\prime_k)\!-\!\!\nabla\! f(\hat{\mathbf{x}}^\prime)) \!-\!2\alpha{\tilde{\mathbf{z}}_{k\!+\!1}}^\intercal\tilde{\mathbf{z}}_k\!+\!\frac{\alpha}{\beta}\Vert{\boldsymbol{\lambda}}_{k\!+\!1}\!-\!{\boldsymbol{\lambda}}_{k}\Vert^2 \notag \\
		\leq&\frac{1}{2} \Vert \Delta{\mathbf{x}}^\prime_{k+1} \Vert^2+(2\alpha^2\ell^2-2\alpha\mu) \Vert \hat{\mathbf{x}}^\prime \Vert^2-\alpha\Vert\tilde{\mathbf{z}}_{k+1}\Vert^2 \notag \\
		+&3\alpha\left(\bar{\sigma}_\mathbf{A}^2\Vert\Delta{\mathbf{x}}_{k+1}\Vert^2\!+\!\bar{\sigma}_\mathcal{L}^2\Vert\Delta{\mathbf{y}}_{k+1}\Vert^2\!+\!\Vert\Delta{\boldsymbol{\delta}}_{k+1}\Vert^2\right) \notag \\
		+&3\alpha\beta(\Vert\tilde{\mathbf{z}}_{k+1}\Vert^2+\eta^2\bar{\sigma}_\mathbf{A}^2\Vert\tilde{\boldsymbol{\lambda}}_{k\!+\!1}\Vert^2+\eta^2\ell^2\Vert\tilde{\mathbf{x}}^\prime_{k+1}\Vert^2), \label{proof_In_DanyRA_7}
	\end{align}
	where we use the Cauchy-Schwarz inequality, the fact that $-2\mathbf{x}^\intercal\mathbf{y}\leq-\Vert\mathbf{x}\Vert^2+\Vert\mathbf{x}-\mathbf{y}\Vert^2$ and $\Vert\mathbf{x}+\mathbf{y}+\mathbf{z}\Vert^2\leq3\Vert\mathbf{x}\Vert^2+3\Vert\mathbf{y}\Vert^2+3\Vert\mathbf{z}\Vert^2$ hold for any vectors $\mathbf{x},\mathbf{y},\mathbf{z}$.
	
	The last line in (\ref{proof_In_DanyRA_6}) can be upper bounded as 
	\begin{align}
		&-\!2\alpha\eta \Vert \mathbf{A}^\intercal \tilde{\boldsymbol{\lambda}}_k \Vert^2\!-2\alpha\eta\tilde{\boldsymbol{\lambda}}_{k}^\intercal\mathbf{A}(\nabla f(\mathbf{x}^\prime_k)-\nabla f(\hat{\mathbf{x}}^\prime)) \notag\\
		\leq&-\alpha\eta\underline{\sigma}_{\mathbf{A}}^2\Vert\tilde{\boldsymbol{\lambda}}_k\Vert^2+\alpha\eta\ell^2\Vert\tilde{\mathbf{x}}^\prime_k\Vert^2. \label{proof_In_DanyRA_9}
	\end{align}
	From the update rule of $\mathbf{x}$, we have $\Vert\mathbf{x}_{k+1}-\mathbf{x}^\prime_{k+1}\Vert^2\leq\Vert\mathbf{x}-\mathbf{x}^\prime_{k+1}\Vert^2$ holds for any $\mathbf{x}$ satisfying the constraint in (\ref{compact_In_DanyRA_5}). Combining this result with the relation $\underline{\sigma}_\mathbf{A}^2\Vert\mathbf{x}-\mathbf{x}^\prime_{k+1}\Vert^2\leq\Vert\mathbf{A}\mathbf{x}-\mathbf{A}\mathbf{x}^\prime_{k+1}\Vert^2$ yields
	\begin{align}
		&\underline{\sigma}_\mathbf{A}^2\Vert\mathbf{x}_{k+1}-\mathbf{x}^\prime_{k+1}\Vert^2 \notag \\
		\leq& \Vert\mathbf{A}\mathbf{x}_{k}\!-\!\gamma(\mathbf{A}\mathbf{x}_k\!+\!\boldsymbol{\mathcal{L}}\mathbf{y}_{k+1}\!+\!\boldsymbol{\delta}_k\!-\!\mathbf{d})\!+\!(\boldsymbol{\delta}_k\!-\!\boldsymbol{\delta}_{k+1})\!-\!\mathbf{A}\mathbf{x}^\prime_{k+1}\Vert^2 \notag \\
		=& \Vert(1\!-\!\gamma)(\tilde{\mathbf{z}}_k\!-\!\tilde{\mathbf{z}}_{k+1}\!+\!\mathbf{A}\mathbf{x}_k\!-\!\mathbf{A}\mathbf{x}^\prime_k\!+\!\mathbf{y}_{k+1}\!-\!\mathbf{y}_{k})\!+\!\gamma\tilde{\mathbf{z}}_{k+1}\Vert^2 \notag \\	\leq&4(1-\gamma)^2\left(3\bar{\sigma}_\mathbf{A}^2\Vert\Delta{\mathbf{x}}_{k+1}\Vert^2\!+\!4\bar{\sigma}_\mathcal{L}^2\Vert\Delta{\mathbf{y}}_{k+1}\Vert^2\!+\!3\Vert\Delta{\boldsymbol{\delta}}_{k+1}\Vert^2\right)\notag \\
		+& 4\left((1-\gamma)^2 \bar{\sigma}_\mathbf{A}^2\Vert\mathbf{x}_k-\mathbf{x}^\prime_k\Vert^2+\gamma^2\Vert\tilde{\mathbf{z}}_{k+1}\Vert^2\right). \label{proof_In_DanyRA_10}
	\end{align}
	Substituting relations (\ref{proof_In_DanyRA_7})-(\ref{proof_In_DanyRA_9}) into (\ref{proof_In_DanyRA_6}) and then combining it with (\ref{proof_In_DanyRA_10})  completes the proof. $\hfill\blacksquare$
	
	\section{Proof of Theorem 1}
	\label{Proof_of_Theorem_1}
	Under conditions in (\ref{conditions_Theorem_1}), it can be verified that all coefficients in the RHS of (\ref{result_Lemma_1}) are negative. For conciseness, we use $-\theta(\theta>0)$ to denote the maximum coefficient among them. In this way, from the relation (\ref{result_Lemma_1}), we have
	\begin{align}
		V_{k+1}
		\leq& V_k\!-\!\theta\left(\Vert\Delta\mathbf{x}^\prime_{k+1}\Vert^2\!\!+\!\Vert\Delta\mathbf{y}_{k+1}\Vert^2 \!\!+\!\Vert\Delta\boldsymbol{\delta}_{k+1}\Vert^2 \!\!+\!\Vert\tilde{\mathbf{x}}^\prime_k\Vert^2\right.\notag \\
		&\left.+\Vert\tilde{\boldsymbol{\lambda}}_k\Vert^2+\Vert\tilde{\mathbf{z}}_k\Vert^2+\Vert\mathbf{x}_k-\mathbf{x}^\prime_k\Vert^2\right). \label{proof_In_DanyRA_11}
	\end{align}
	It can be derived from (\ref{proof_In_DanyRA_11}) and the relation $\hat{\mathbf{x}}=\hat{\mathbf{x}}^\prime$ that
	\begin{align*}
		\sum_{k=0}^\infty\Vert\mathbf{x}_k\!-\!\hat{\mathbf{x}}\Vert^2\!\leq\!\sum_{k=0}^\infty2\left(\Vert\mathbf{x}_k\!-\!\mathbf{x}_k^\prime\Vert^2+\Vert\tilde{\mathbf{x}}^\prime_k\Vert^2\right)\leq\frac{2V_0}{\theta}<\infty,
	\end{align*}
	which indicates that $\lim_{k\rightarrow\infty}\Vert\mathbf{x}_k-\hat{\mathbf{x}}\Vert=0$. Combining this relation with (\ref{Sensitivity_result}) proves  the result \romannumeral1).
	
	Then we will prove the feasibility-related results \romannumeral2) and \romannumeral3). According to the definition of the feasible set $\mathcal{K}^x(\mathbf{x}_{k})$ in (\ref{compact_In_DanyRA_5}), we have 
	\begin{align*}
		\mathbf{A}\mathbf{x}_{k+1}\!=\!\mathbf{A}\mathbf{x}_{k}\!-\!\gamma(\mathbf{A}\mathbf{x}_k\!+\!\boldsymbol{\mathcal{L}}\mathbf{y}_{k+1}\!+\!\boldsymbol{\delta}_k\!-\!\mathbf{d})\!+\!(\boldsymbol{\delta}_k\!-\!\boldsymbol{\delta}_{k+1}). 
	\end{align*}
	Subtracting $\mathbf{d}$ from both sides of the above relation and then multiplying by 
	$\hat{\mathbf{1}}^\intercal$, we have that
	\begin{align*}
		\sum_{i=1}^n \left(A_ix_{i,k+1}\!+\!\delta_{i,k+1}\!-\!d_{i}\right)\!=\!(1\!-\!\gamma)\sum_{i=1}^n \left(A_ix_{i,k}\!+\!\delta_{i,k}\!-\!d_{i}\right)
	\end{align*}
	holds for any $k\geq 0$. {If  $\sum_{i=1}^n \left(A_ix_{i,k}+\delta_{i,k}-d_{i}\right)\preccurlyeq \frac{n\omega \mathbf{1}}{1-\gamma}$ holds at the $k_0$-th iteration, we have}
	{
		\begin{align}
			\sum_{i=1}^n \left(A_ix_{i,k_0+1}-d_{i}\right)\preccurlyeq  n\omega\mathbf{1}-\sum_{i=1}^n\delta_{i,k_0+1}\preccurlyeq 0, \label{proof_In_DanyRA_16}
		\end{align}
		which can be extended to any $k\!\geq\! k_0$, proving \romannumeral3).
		
		On the other hand, when there exists a large violation  at the $k_0$-th iteration, it can be derived that
		\begin{align}
			&\sum_{i=1}^n \left(A_ix_{i,k_0+t}+\delta_{i,k_0+t}-d_{i}\right) \notag \\ \preccurlyeq & (1-\gamma)^t (\sum_{i=1}^n  \left(A_ix_{i,k_0}-d_{i}\right)+\sum_{i=1}^n\delta_{i,k_0}) \notag \\
			\preccurlyeq & (1-\gamma)^t(C_{k_0}^{\rm vio}+C_\delta)\mathbf{1} \label{proof_In_DanyRA_18}
		\end{align}	
		It can be derived from \eqref{proof_In_DanyRA_11} that $\|\boldsymbol\delta_k-\hat{\boldsymbol\delta}\|^2
		\leq V_k\leq V_0$. Combining this result with the fact that   $\left\|\sum_{i=1}^n\delta_{i,k_0}\right\|_\infty
		\leq\sqrt n\|\boldsymbol\delta_{k_0}\|$ yields
		\begin{align}
			\left\|\sum_{i=1}^n\delta_{i,k_0}\right\|_\infty
			\leq\sqrt n(\|\hat{\boldsymbol\delta}\|+\|\boldsymbol\delta_{k_0}-\hat{\boldsymbol\delta}\|)\leq C_\delta. \label{proof_In_DanyRA_19}
		\end{align}	
		Combining relations \eqref{proof_In_DanyRA_18} with \eqref{proof_In_DanyRA_19} yields
		\begin{align}
			\sum_{i=1}^n \! \left(A_ix_{i,k_0+t}\!+\!\delta_{i,k_0+t}\!-\!d_{i}\right)
			\!\preccurlyeq\!  (1\!-\!\gamma)^t(C_{k_0}^{\rm vio}\!+\!C_\delta)\mathbf{1}. \label{proof_In_DanyRA_20}
	\end{align}}	
	From \romannumeral3), when the RHS of (\ref{proof_In_DanyRA_20}) is less than or equal to $\frac{n\omega \mathbf{1}}{1-\gamma}$, the global constraint is satisfied thereafter, and thus the result \romannumeral2) is proved.
	The proof is completed. $\hfill\blacksquare$
	
	\section{Proof of Theorem 2}
	\label{Proof_of_Theorem_2}
	In the case involving a virtual queue with decaying buffer $\omega_{k}$, the $\boldsymbol{\delta}$-update rule can be rewritten as 
	\begin{align} 	
		\boldsymbol{\delta}_{k+1}=\max\{\boldsymbol{\delta}_k-\alpha(\mathbf{z}_k-\mathbf{d}+\boldsymbol{\lambda}_k),{\omega_k}\}, \label{compact_In_DanyRA_3_k}
	\end{align}
	Similar to the proof of Lemma 1, from (\ref{compact_In_DanyRA_3_k}), we have
	\begin{align}
		&\Vert\tilde{\boldsymbol{\delta}}_{k+1}\Vert^2-\Vert\tilde{\boldsymbol{\delta}}_{k}\Vert^2 \notag\\ 
		=&-\Vert\Delta{\boldsymbol{\delta}}_{k+1}\Vert^2-2\alpha \tilde{\boldsymbol{\delta}}_{k+1}^\intercal (\mathbf{z}_k+\boldsymbol{\lambda}_k-\mathbf{d})\notag \\
		&+\!2\left(\tilde{\boldsymbol{\delta}}_{k+1}\!+\!\omega_k\mathbf{1}\!-\!\omega_k\mathbf{1}\right)^\intercal\left(\boldsymbol{\delta}_{k+1}\!-\!(\boldsymbol{\delta}_{k}\!-\!\alpha(\mathbf{z}_k\!-\!\mathbf{d}\!+\!\boldsymbol{\lambda}_k))\right) \notag \\
		\leq & -\Vert\Delta{\boldsymbol{\delta}}_{k+1}\Vert^2+2(\omega_k\mathbf{1})^\intercal\left(\boldsymbol{\delta}_{k+1}-\boldsymbol{\delta}_{k}+\alpha(\tilde{\mathbf{z}}_k+\tilde{\boldsymbol{\lambda}}_k)\right)\notag \\
		&-2\alpha \tilde{\boldsymbol{\delta}}_{k+1}^\intercal (\tilde{\mathbf{z}}_k+\tilde{\boldsymbol{\lambda}}_k) -2\alpha(\tilde{\boldsymbol{\delta}}_{k+1}-\omega_k\mathbf{1})^\intercal\hat{\boldsymbol{\lambda}}. \label{proof_In_DanyRA_4_k}
	\end{align}
	
	The last term in the RHS of (\ref{proof_In_DanyRA_4_k}) can be rewritten as 
	\begin{align}
		-2\alpha({\boldsymbol{\delta}}_{k+1}-\omega_k\mathbf{1})^\intercal\hat{\boldsymbol{\lambda}}+2\alpha\hat{\boldsymbol{\delta}}^\intercal\hat{\boldsymbol{\lambda}}
		\leq 0, \label{proof_In_DanyRA_5_k}
	\end{align}
	where we use ${\boldsymbol{\delta}}_{k+1}-\omega_k\mathbf{1} \succcurlyeq 0$, $\hat{\boldsymbol{\lambda}}\succcurlyeq 0$, $\hat{\boldsymbol{\delta}}\succcurlyeq 0$, and ${\hat{\boldsymbol{\lambda}}}^\intercal\hat{\boldsymbol{\delta}}=0$. 
	
	Similar to the proof of Lemma 1, combining relations (\ref{proof_In_DanyRA_1})-(\ref{proof_In_DanyRA_3}), (\ref{proof_In_DanyRA_4_k}) and (\ref{proof_In_DanyRA_5_k}) yields
	\begin{align}
		&V_{k+1}^{\prime}-V_k^{\prime} 
		\leq  -\Vert\Delta{\mathbf{x}}^\prime_{k+1}\Vert^2-\Vert\Delta{\mathbf{y}}_{k+1}\Vert^2-\Vert\Delta{\boldsymbol{\delta}}_{k+1}\Vert^2 \notag\\
		&  -2\alpha\tilde{\mathbf{x}}_{k+1}^{\prime\intercal}(\nabla f(\mathbf{x}_k^\prime)\!-\!\nabla f(\hat{\mathbf{x}}^\prime))\!-\!2\alpha{\tilde{\mathbf{z}}_{k+1}}^\intercal\tilde{\mathbf{z}}_k  \!+\!\frac{\alpha}{\beta}\Vert{\boldsymbol{\lambda}}_{k+1}\!-\!{\boldsymbol{\lambda}}_{k}\Vert^2 \notag\\
		&-2\alpha\eta \Vert \mathbf{A}^\intercal \tilde{\boldsymbol{\lambda}}_k \Vert^2-2\alpha\eta\tilde{\boldsymbol{\lambda}}_{k}^\intercal\mathbf{A}(\nabla f(\mathbf{x}^\prime_k)-\nabla f(\hat{\mathbf{x}}^\prime)) \notag \\
		&+2(\omega_k\mathbf{1})^\intercal\left(\boldsymbol{\delta}_{k+1}-\boldsymbol{\delta}_{k}+\alpha(\tilde{\mathbf{z}}_k+\tilde{\boldsymbol{\lambda}}_k)\right). \label{proof_In_DanyRA_6_k}
	\end{align}
	Then, we will give the upper bound of the RHS of (\ref{proof_In_DanyRA_6_k}) to prove the convergence of the algorithm. The last line in (\ref{proof_In_DanyRA_6_k}) can be upper bounded as 
	\begin{align}
		&2(\omega_k\mathbf{1})^\intercal\left(\boldsymbol{\delta}_{k+1}-\boldsymbol{\delta}_{k}+\alpha(\tilde{\mathbf{z}}_k+\tilde{\boldsymbol{\lambda}}_k)\right) \notag \\
		\leq&\frac{1}{2}\left(\Vert\Delta\boldsymbol{\delta}_{k+1}\Vert^2+\alpha^2\Vert\tilde{\mathbf{z}}_{k}\Vert^2+\alpha^2\Vert\tilde{\boldsymbol{\lambda}}_{k}\Vert^2\right)+6n\omega_k^2. \label{proof_In_DanyRA_8_k}
	\end{align}
	
	Substituting relations (\ref{proof_In_DanyRA_7}), (\ref{proof_In_DanyRA_9}), (\ref{proof_In_DanyRA_8_k}) into (\ref{proof_In_DanyRA_6_k}) and then combining it with (\ref{proof_In_DanyRA_10}) yield 
	\begin{align}
		&V_{k+1}-V_k-6n\omega_k^2	\notag\\
		\leq& \frac{6\alpha\bar{\sigma}_\mathbf{A}^2(1\!+\!3c)\!-\!1}{2}\Vert\Delta\mathbf{x}^\prime_{k\!+\!1}\Vert^2\!+\!(3\alpha\bar{\sigma}_\mathcal{L}^2(1\!+\!4c)\!-\!1)\Vert\Delta\mathbf{y}_{k\!+\!1}\Vert^2 \notag \\
		+&\!\frac{6\alpha(1\!+\!3c)\!-\!1}{2}\Vert\Delta\boldsymbol{\delta}_{k\!+\!1}\Vert^2 \!\!+\! \alpha(2\ell^2\alpha\!-\!2\mu\!+\!\eta\ell^2(3\beta\eta\!+\!1))\Vert\tilde{\mathbf{x}}^\prime_k\Vert^2\notag \\
		+&\alpha\big(\eta(3\beta\eta\bar{\sigma}_\mathbf{A}^2-\underline{\sigma}_\mathbf{A}^2)+\frac{\alpha}{2}\big)\Vert\tilde{\boldsymbol{\lambda}}_k\Vert^2+\frac{\alpha(3\beta+\alpha-1)}{2}\Vert\tilde{\mathbf{z}}_k\Vert^2\notag \\
		+&\frac{\underline{\sigma}_\mathbf{A}^2\alpha(1\!-\!3\beta)}{8\gamma^2}\left({4(1-\gamma)^2\kappa_\mathbf{A}^2}\!-\!1\right)\Vert\mathbf{x}_k-\mathbf{x}^\prime_k\Vert^2. \label{proof_In_DanyRA_11_k}
	\end{align}
	Under conditions in (\ref{conditions_Theorem_1}), it can be verified that all coefficients in the RHS of (\ref{result_Lemma_1}) are negative. For conciseness, we use $-\theta(\theta>0)$ to denote the maximum coefficient among them. In this way, from the relation (\ref{result_Lemma_1}), we have
	\begin{align}
		V_{k+1}
		&\leq V_k-\theta\left(\Vert\Delta\mathbf{x}^\prime_{k+1}\Vert^2\!+\!\Vert\Delta\mathbf{y}_{k+1}\Vert^2 +\Vert\Delta\boldsymbol{\delta}_{k+1}\Vert^2 \!\!\right.\notag \\
		&\left.+\Vert\tilde{\mathbf{x}}^\prime_k\Vert^2\!\!+\!\Vert\tilde{\boldsymbol{\lambda}}_k\Vert^2\!\!+\!\Vert\tilde{\mathbf{z}}_k\Vert^2\!\!+\!\Vert\mathbf{x}_k\!-\!\mathbf{x}^\prime_k\Vert^2\!\right)\!\!+\!6n\omega_k^2. \label{proof_In_DanyRA_12_k}
	\end{align}
	Combining (\ref{proof_In_DanyRA_12_k}) with $\hat{\mathbf{x}}=\hat{\mathbf{x}}^\prime$, and $\sum_{k=0}^\infty\omega_k^2<\infty$ yields
	\begin{align*}
		\sum_{k=0}^\infty\Vert\mathbf{x}_k-\hat{\mathbf{x}}\Vert^2&\leq\sum_{k=0}^\infty2\left(\Vert\mathbf{x}_k-\mathbf{x}^\prime_k\Vert^2+\Vert\tilde{\mathbf{x}}^\prime_k\Vert^2\right)\notag\\
		&\leq2V_0+12n\sum_{k=0}^\infty\omega_k^2<\infty,
	\end{align*}
	which indicates that $\lim_{k\rightarrow\infty}\Vert\mathbf{x}_k-\hat{\mathbf{x}}\Vert=0$. Combining this relation with (\ref{Sensitivity_result}) proves  the result \romannumeral1).
	
	The proofs for feasibility-related results \romannumeral2) and \romannumeral3) follow analogously to Theorem 1, and thus the detailed derivations are omitted. The proof is completed. $\hfill\blacksquare$
	
	\section{Proof of Theorem 3}
	\label{Proof_of_Theorem_3}
	Since Eq-DanyRA algorithm can be regarded as a special case of the DanyRA algorithm  with known $\boldsymbol{\delta}^\star$, it can be derived from Lemma 1 that
	\begin{align}
		&V_{k+1}-V_k	\notag\\
		&\leq \frac{6\alpha\bar{\sigma}_\mathbf{A}^2(1\!+\!3c)\!-\!1}{2}\Vert\Delta\mathbf{x}^\prime_{k\!+\!1}\Vert^2\!+\!(3\alpha\bar{\sigma}_\mathcal{L}^2(1\!+\!4c)\!-\!1)\Vert\Delta\mathbf{y}_{k\!+\!1}\Vert^2 \notag \\
		&+\! \alpha(2\ell^2\alpha\!-\!2\mu\!+\!\eta\ell^2(3\beta\eta\!+\!1))\Vert\tilde{\mathbf{x}}^\prime_k\Vert^2+\frac{\alpha(3\beta-1)}{4}\Vert\tilde{\mathbf{z}}_k\Vert^2\notag \\
		&+\alpha\eta(3\beta\eta\bar{\sigma}_\mathbf{A}^2-\underline{\sigma}_\mathbf{A}^2)\Vert\tilde{\boldsymbol{\lambda}}_k\Vert^2+\frac{\alpha(3\beta-1)}{4}\Vert\tilde{\mathbf{z}}_k\Vert^2\notag \\
		&+\frac{\underline{\sigma}_\mathbf{A}^2\alpha(1\!-\!3\beta)}{8\gamma^2}\left({4(1-\gamma)^2\kappa_\mathbf{A}^2}\!-\!1\right)\Vert\mathbf{x}_k-\mathbf{x}^\prime_k\Vert^2. \label{Lemma_2_1}
	\end{align}
	Under conditions in (\ref{conditions_Theorem_1}), we have $\frac{\alpha(3\beta-1)}{2}<0$. Then, according to the definition of $\mathbf{z}$, we have 
	\begin{align}
		\frac{\alpha(3\beta\!\!-\!\!1)}{4}\Vert\tilde{\mathbf{z}}_k\Vert^2
		\!\!\leq\!\! \frac{\alpha(1\!-\!3\beta)}{4}\Big(\bar{\sigma}^2_\mathbf{A}\Vert\tilde{\mathbf{x}}^\prime_k\Vert^2\!\!-\!\!\frac{\underline{\sigma}_\mathcal{L}^2}{2}\Vert\tilde{\mathbf{y}}_k\Vert^2\Big), \label{Lemma_2_2}
	\end{align}
	where the Cauchy-Schwarz inequality is used.
	
	Substituting (\ref{Lemma_2_2}) into the second line of (\ref{Lemma_2_1}), it can be verified that all coefficients in the RHS of the derived relation are negative under conditions in (\ref{conditions_Theorem_1}) and (\ref{conditions_Lemma_2}). Then, arranging the formulas yields
	\begin{align}
		V_{k+1}	\leq \theta^\prime V_k. \label{result_Lemma_2}
	\end{align}
	
	It can be derived from  the definition of $V_k$ that
	\begin{align*}
		\Vert\mathbf{x}_k\!-\!\hat{\mathbf{x}}\Vert^2\!\leq\!2\left(\Vert\tilde{\mathbf{x}}^\prime_k\Vert^2\!+\!\Vert\mathbf{x}_k\!-\!\mathbf{x}^\prime_k\Vert^2\right) \leq\frac{2V_k}{\min\{1,\frac{\underline{\sigma}_\mathbf{A}^2\alpha(1-3\beta)}{8\gamma^2}\}}.
	\end{align*}
	Combining this relation with (\ref{result_Lemma_2}) proves the result \romannumeral1).

	Next, subtracting $d_i$ from both sides of the equality constraint in the $\mathbf{x}$-update step of Eq-DanyRA algorithm and then summing the result from $i=1$ to $n$ yields
	\begin{align*}
		\sum_{i=1}^n \left(A_ix_{i,k+1}\!-\!d_{i}\right)\!=\!(1\!-\!\gamma)\sum_{i=1}^n \left(A_ix_{i,k}\!-\!d_{i}\right),  \forall k\geq 0,
	\end{align*}
	which proves  \romannumeral2) and \romannumeral3). The proof is completed. $\hfill\blacksquare$
	
	\bibliographystyle{plain}        
	\bibliography{autosam}           

\end{document}